# Analytic properties of zeta functions and subgroup growth

By Marcus du Sautoy and Fritz Grunewald

## 1. Introduction

It has become somewhat of a cottage industry over the last fifteen years to understand the rate of growth of the number of subgroups of finite index in a group $G$. Although the story began much before, the recent activity grew out of a paper by Dan Segal in [36]. The story so far has been well-documented in Lubotzky's subsequent survey paper in [30].

In [24] the second author of this article, Segal and Smith introduced the zeta function of a group as a tool for understanding this growth of subgroups. Let $a_n(G)$ be the number of subgroups of index $n$ in the finitely generated group $G$ and $s_N(G) = a_1(G) + \cdots + a_N(G)$ be the number of subgroups of index $N$ or less. The zeta function is defined as the Dirichlet series with coefficients $a_n(G)$ and has a natural interpretation as a noncommutative generalization of the Dedekind zeta function of a number field:

$$
\begin{aligned}
(1.1) \qquad \zeta_G(s) &= \sum_{n=1}^{\infty} a_n(G) n^{-s} \\
&= \sum_{H \leq G} |G : H|^{-s}.
\end{aligned}
$$

For example, without such a tool it would be difficult to prove that the number of subgroups in the rank-two free abelian group $G = \mathbb{Z}^2$ grows as follows:

$$
s_N(\mathbb{Z}^2) \sim \left( \pi^2/12 \right) N^2
$$

as $N$ tends to infinity. (Here $f(n) \sim g(n)$ means $f(n)/g(n)$ tends to 1 as $n$ tends to infinity.) This is a consequence of the expression for the zeta function of the free abelian group of rank $d$:

$$
\zeta_{\mathbb{Z}^d}(s) = \zeta(s) \cdots \zeta(s - d + 1)
$$

where $\zeta(s)$ is the Riemann zeta function.

The zeta function (1.1) defines an analytic function on some right half of the complex plane $\Re(s) > \alpha(G)$ precisely when the coefficients $a_n(G)$ are bounded by a polynomial. A characterization of such finitely generated residually finite groups, groups of polynomial subgroup growth, was provided by Lubotzky, Mann and Segal [31]. They are groups which have a subgroup of finite index that is soluble of finite rank.

In this paper we consider the analytic behaviour of the zeta function of the subclass of finitely generated nilpotent groups. This class of groups has the added bonus that their zeta functions satisfy an Euler product (see [24]):

$$\zeta_G(s) = \prod_p \zeta_{G,p}(s)$$

where the local factors for each prime $p$ are defined as:

$$\zeta_{G,p}(s) = \sum_{n=0}^{\infty} a_{p^n}(G) p^{-ns}.$$

It was also proved in [24] that if the nilpotent group is torsion-free then these local factors are all rational functions in $p^{-s}$. However the proof gave little understanding of how these rational functions varied as $p$ varied and was not sufficient to understand the global behaviour of $\zeta_G(s)$.

In this paper we introduce some new methods to understand the analytic behaviour of the zeta function of a group. We can then combine this knowledge with suitable Tauberian theorems to deduce results about the growth of subgroups in a nilpotent group. In order to state our results we introduce the following notation. For $\alpha \in \mathbb{R}$ and $N \in \mathbb{N}$, define

$$s_N^\alpha(G) := \sum_{n=1}^{N} \frac{a_n(G)}{n^\alpha}.$$

We prove the following:

THEOREM 1.1. *Let $G$ be a finitely generated nilpotent infinite group.*

(1) *The abscissa of convergence $\alpha(G)$ of $\zeta_G(s)$ is a rational number and $\zeta_G(s)$ can be meromorphically continued to $\Re(s) > \alpha(G) - \delta$ for some $\delta > 0$. The continued function is holomorphic on the line $\Re(s) = (\alpha)G$ except for a pole at $s = \alpha(G)$.*

(2) *There exist a nonnegative integer $b(G) \in \mathbb{N}$ and some real numbers $c, c' \in \mathbb{R}$ such that*

$$s_N(G) \quad \sim \quad c \cdot N^{\alpha(G)} (\log N)^{b(G)}$$
$$s_N^{\alpha(G)}(G) \quad \sim \quad c' \cdot (\log N)^{b(G)+1}$$

*for $N \to \infty$.*

Whether the abscissa of convergence is a rational number was raised as one of the major open problems in the field in Lubotzky's survey article [30].

Note that the integer $b(G) + 1$ is the multiplicity of the pole of $\zeta_G(s)$ at $s = \alpha(G)$. In [13] several examples are given where this multiplicity is greater than one. For example, the zeta function of the discrete Heisenberg group

$$G = \begin{pmatrix} 1 & \mathbb{Z} & \mathbb{Z} \\ 0 & 1 & \mathbb{Z} \\ 0 & 0 & 1 \end{pmatrix}$$

has the following expression:

$$\zeta_G(s) = \zeta(s)\zeta(s-1)\zeta(2s-2)\zeta(2s-3) \cdot \zeta(3s-3)^{-1}.$$

The double pole at $s = 2$ implies that the growth of subgroups is:

$$s_N(G) \sim \frac{\zeta(2)^2}{2\zeta(3)} N^2 \log N$$

for $N \to \infty$. This was first observed in Smith's thesis [37]. This example has meromorphic continuation to the whole complex plane. In [22] it is shown that this is also true for any finite extension of a free abelian group. In general, though, these functions have natural boundaries as discussed in [13]. However we have introduced in a separate paper [16] the concept of the ghost zeta function which does tend to have meromorphic continuation.

The proof of the meromorphic continuation of the zeta function of a nilpotent group depends on showing a more general result which holds for any zeta function which can be defined as an Euler product over primes $p$ of *cone integrals over* $\mathbb{Q}$.

*Definition* 1.2. (1) Let $\psi(\mathbf{x})$ be a formula in the first order language (in the sense of logic) for the valued field $\mathbb{Q}_p$ built from the following symbols: $+$ (addition), $\cdot$ (multiplication), $|$ (here $x|y$ means $v(x) \leq v(y)$), for every element of $\mathbb{Q}_p$ a symbol denoting that element, $=$, $\wedge$ (and), $\vee$ (or), $\neg$ (not), and quantifiers $\exists x$ (there exists $x \in \mathbb{Q}_p$ :) and $\forall x$ (for every $x \in \mathbb{Q}_p$ :).

The formula $\psi(\mathbf{x})$ is called *a cone condition over* $\mathbb{Q}$ if there exist nonzero polynomials $f_i(\mathbf{x}), g_i(\mathbf{x})(i = 1, \ldots, l)$ over $\mathbb{Q}$ in the variables $\mathbf{x} = x_1, \ldots, x_m$ such that $\psi(\mathbf{x})$ is a conjunction of formulas

$$v(f_i(\mathbf{x})) \leq v(g_i(\mathbf{x}))$$

for $i = 1, \ldots, l$.

(2) Given a cone condition $\psi(\mathbf{x})$ over $\mathbb{Q}$ and nonzero polynomials $f_0$ and $g_0$ with coefficients in $\mathbb{Q}$, we call an integral

$$Z_{\mathcal{D}}(s, p) = \int_{V_p = \{\mathbf{x} \in \mathbb{Z}_p^m : \psi(\mathbf{x}) \text{ is valid}\}} |f_0(\mathbf{x})|^s |g_0(\mathbf{x})| \, |dx|$$

*a cone integral defined over* $\mathbb{Q}$, where $|dx|$ is the normalized additive Haar measure on $\mathbb{Z}_p^m$ and $\mathcal{D} = \{f_0, g_0, f_1, g_1, \ldots, f_l, g_l\}$ is called the *cone integral data*.

(3) We say that a function $Z(s)$ is *defined as an Euler product of cone integrals over* $\mathbb{Q}$ with cone integral data $\mathcal{D}$ if

$$Z(s) = Z_{\mathcal{D}}(s) = \prod_{p \text{ prime, } a_{p,0} \neq 0} \left( a_{p,0}^{-1} \cdot Z_{\mathcal{D}}(s, p) \right)$$

where $a_{p,0} = Z_{\mathcal{D}}(\infty, p)$ is the constant coefficient of $Z_{\mathcal{D}}(s, p)$; i.e., we normalize the local factors to have constant coefficient 1.

We shall explain during the course of the analysis of cone integrals why $a_{p,0} \neq 0$ for almost all primes $p$.

In Section 5 we show that for a nilpotent group $G$, $\zeta_G(s) = Z_{\mathcal{D}}(s-d) \cdot P(s)$ where $Z_{\mathcal{D}}(s)$ is defined as an Euler product of cone integrals over $\mathbb{Q}$, $P(s) = \prod_{p \in S} P_p(p^{-s})$ where $S$ is a finite set of primes, $P_p(X)$ is a rational function and $d$ is the Hirsch length of $G$. (The Hirsch length is the number of infinite cyclic factors in a composition series for $G$.)

We adapt some ideas of Denef introduced in [5] to give an explicit expression for a cone integral, valid for almost all primes $p$ in terms of the resolution of singularities $(Y, h)$ of the polynomial $F(\mathbf{x}) = \prod_{i=0}^{l} f_i(\mathbf{x})g_i(\mathbf{x})$. In particular we show:

THEOREM 1.3. *Let* $(Y, h)$ *be a resolution over* $\mathbb{Q}$ *for* $F(\mathbf{x}) = \prod_{i=0}^{l} f_i(\mathbf{x})g_i(\mathbf{x})$ *and let* $E_i, i \in T$, *be the irreducible components of the reduced scheme* $(h^{-1}(D))_{\text{red}}$ *over* $\mathrm{Spec}(\mathbb{Q})$ *where* $D = \mathrm{Spec}\left(\frac{\mathbb{Q}[\mathbf{x}]}{(F)}\right)$. *Then there exist rational functions* $P_I(x, y) \in \mathbb{Q}(x, y)$ *for each* $I \subset T$ *with the property that for almost all primes* $p$

$$(1.2) \qquad\qquad Z_{\mathcal{D}}(s, p) = \sum_{I \subset T} c_{p,I} P_I(p, p^{-s})$$

*where*

$$c_{p,I} = \mathrm{card}\{a \in \overline{Y}(\mathbb{F}_p) : a \in \overline{E_i} \text{ if and only if } i \in I\}$$

*and* $\overline{Y}$ *means the reduction* mod $p$ *of the scheme* $Y$.

The $E_i$ are smooth quasiprojective varieties defined over $\mathbb{Q}$ and we can use the Lang-Weil estimates for the number of points on such varieties mod $p$ to identify the abscissa of convergence of the global zeta function $Z_{\mathcal{D}}(s)$.

However just knowing the shape of the zeta function from the expression (1.2) is not sufficient to infer that the Euler product of these expressions can be meromorphically continued beyond its region of convergence. For example,

$$(1.3) \qquad\qquad \prod_{p \text{ prime}} \left( 1 + \frac{p^{-1-s}}{(1 - p^{-s})} \right)$$

converges for $\Re(s) > 0$ but has $\Re(s) = 0$ as a natural boundary. We give instead a subtler expression for the cone integrals. Rather than a sum over the subsets of $T$, the indexing set of the irreducible components, this second expression is a sum over the open simplicial pieces of a natural polyhedral cone that one associates to the cone condition $\psi$.

THEOREM 1.4. *There exist a closed polyhedral cone $\overline{D}$ in $\mathbb{R}_{\geq 0}^t$ where $t = \operatorname{card} T$ and a simplicial decomposition into open simplicial pieces denoted by $R_k$ where $k \in \{0, 1, \ldots, w\}$. Let $R_0 = (0, \ldots, 0)$ and $R_1, \ldots, R_q$ be the one-dimensional pieces. For each $k \in \{0, 1, \ldots, w\}$ let $M_k \subset \{1, \ldots, q\}$ denote those one-dimensional pieces in the closure $\overline{R_k}$ of $R_k$. Then there exist positive integers $A_j, B_j$ for $j \in \{1, \ldots, q\}$ such that for almost all primes $p$*

$$(1.4) \qquad Z_{\mathcal{D}}(s, p) = \sum_{k=0}^{w} (p-1)^{I_k} p^{-m} c_{p, I_k} \prod_{j \in M_k} \frac{p^{-(A_j s + B_j)}}{\left(1 - p^{-(A_j s + B_j)}\right)}$$

*where $c_{p, I_k}$ is as defined in Theorem 1.3 and $I_k$ is the subset of $T$ defined so that $i \in T \backslash I_k$ if and only if the $i^{\text{th}}$ coordinate is zero for all elements of $R_k$.*

This expression (1.4) motivates the name *cone integral*. An explicit expression is given for the integers $A_j$ and $B_j$ in terms of the numerical data of the resolution. It is contained in the proof of this theorem which appears in Section 3. At the end of Section 3 we also give an expression for the rational functions of cone integrals at primes with bad reduction, which shows that they are not far from the expression in (1.4). In particular, the local poles at bad primes are a subset of the candidate poles $-B_j / A_j$, $j \in \{1, \ldots, q\}$ provided by the expression (1.4) for good primes.

With this combinatorial expression in hand, we can show that the pathologies of examples like (1.3) do not arise. In particular, we show that the abscissa of convergence of the global zeta function is determined by the terms in the expression (1.4) corresponding to the one-dimensional edges $R_1, \ldots, R_q$. We then show how to use Artin $L$-functions to analytically continue a function like

$$\prod_{p \text{ prime}} \left(1 + c_{p, I_k} \frac{p^{-s}}{(1 - p^{-s})}\right)$$

beyond its region of convergence. We can then use various Tauberian theorems to estimate the growth of the coefficients in the Dirichlet series expressing $Z_{\mathcal{D}}(s)$. In particular we prove the following:

THEOREM 1.5. *Let $Z(s)$ be defined as an Euler product of cone integrals over $\mathbb{Q}$. Then $Z(s)$ is expressible as a Dirichlet series $\sum_{n=1}^{\infty} a_n n^{-s}$ with nonnegative coefficients $a_n$. Suppose that $Z(s)$ is not the constant function.*

(1) *The abscissa of convergence $\alpha$ of $Z(s)$ is a rational number and $Z(s)$ has a meromorphic continuation to $\Re(s) > \alpha - \delta$ for some $\delta > 0$. The continued function is holomorphic on the line $\Re(s) = \alpha$ except for a pole at $s = \alpha$.*

(2) *Let the pole at $s = \alpha$ have order $w$. Then there exist some real numbers $c, c' \in \mathbb{R}$ such that*

$$a_1 + a_2 + \cdots + a_N \quad \sim \quad c \cdot N^\alpha (\log N)^{w-1}$$
$$a_1 + a_2 2^{-\alpha} + \cdots + a_N N^{-\alpha} \quad \sim \quad c' \cdot (\log N)^w$$

*for $N \to \infty$.*

One of the key problems in this area was to link zeta functions of groups up to questions in some branch of more classical number theory. We restate Theorem 1.3 explicitly for groups as it provides just such a path from zeta functions of groups to the more classical question of counting points mod $p$ on a variety. The path is quite explicit. We define in Section 5 a polynomial $F_G$ over $\mathbb{Q}$ associated to each nilpotent group $G$.

THEOREM 1.6. *Let $G$ be a finitely generated nilpotent group. Let $(Y, h)$ be a resolution over $\mathbb{Q}$ for the polynomial $F_G$. Let $E_i, i \in T$ be the irreducible components of the reduced scheme $(h^{-1}(D))_{\mathrm{red}}$ associated to $h^{-1}(D)$ where $D = \mathrm{Spec}\left(\frac{\mathbb{Q}[\mathbf{x}]}{(F_G)}\right)$. Then there exist rational functions $P_I(x, y) \in \mathbb{Q}(x, y)$ for each $I \subset T$ with the property that for almost all primes $p$*

$$\zeta_{G,p}(s) = \sum_{I \subset T} c_{p,I} P_I(p, p^{-s})$$

*where*

$$c_{p,I} = \mathrm{card}\{a \in \overline{Y}(\mathbb{F}_p) : a \in \overline{E_i} \text{ if and only if } i \in I\}$$

*and $\overline{Y}$ means the reduction mod $p$ of the scheme $Y$.*

The behaviour of the local factors as we vary $p$ is one of the other major problems in the field. For example in the Heisenberg group with entries from a quadratic number field, the behaviour of the local factors depends on how $p$ behaves in the number field [24]. Our explicit formula however takes the subject away from the behaviour of primes in number fields to the problem of counting points mod $p$ on a variety, a question which is in general wild and far from the uniformity predicted by all previous examples (see [24] and [15]). Two papers [11] and [12] by the first author contain an example of a class two nilpotent group of Hirsch length 9 whose zeta function depends on counting points mod $p$ on the elliptic curve $y^2 = x^3 - x$.

Given a nilpotent group $G$ it is possible to construct and analyse the polynomial $F_G$ in question. For example, in the free abelian group or the Heisenberg group, the polynomial does not require any resolution of singularities, as $D$ in this case only involves normal crossings. Hence the $E_i, i \in T$, in this case are just the irreducible components of the algebraic set $F_G = 0$. However this is not true in general. For example the class two nilpotent group defined using the elliptic curve mentioned above has an $F_G$ whose singularities are not normal crossings and which therefore require some resolution.

We have put the emphasis in this introduction on applying these cone integrals to the question of counting subgroups in nilpotent groups; however our results extend in a number of other directions.

(1) Variants of our zeta functions have been considered which count only subgroups with some added feature, for example normal subgroups. Define

$$
\begin{aligned}
a_n^{\triangleleft}(G) &= \operatorname{card}\{H : H \text{ is normal subgroup of } G \text{ and } |G : H| = n\}, \\
\zeta_G^{\triangleleft}(s) &= \sum a_n^{\triangleleft}(G) n^{-s}.
\end{aligned}
$$

Our theorems hold for this normal zeta function and many of the other variants.

(2) Let $L$ be a ring additively isomorphic to $\mathbb{Z}^d$. Define

$$
\begin{aligned}
a_n(L) &= \operatorname{card}\{H : H \text{ is a subring of } L \text{ and } |L : H| = n\}, \\
a_n^{\triangleleft}(L) &= \operatorname{card}\{H : H \text{ is an ideal of } L \text{ and } |L : H| = n\}.
\end{aligned}
$$

Zeta functions of $L$ were also defined in [24] as the Dirichlet series

$$
\begin{aligned}
\zeta_L(s) &= \sum a_n(L) n^{-s}, \\
\zeta_L^{\triangleleft}(s) &= \sum a_n^{\triangleleft}(L) n^{-s}.
\end{aligned}
$$

It was pointed out in [24] that these zeta functions have an Euler product; as for the case of nilpotent groups:

$$
\begin{aligned}
\zeta_L(s) &= \prod_{p \text{ prime}} \zeta_{L \otimes \mathbb{Z}_p}(s), \\
\zeta_L^{\triangleleft}(s) &= \prod_{p \text{ prime}} \zeta_{L \otimes \mathbb{Z}_p}^{\triangleleft}(s).
\end{aligned}
$$

Unlike the situation for groups, there is no need to make an assumption of nilpotency in the case of rings. We can therefore consider examples like $L = \mathfrak{sl}_2(\mathbb{Z})$ or the $\mathbb{Z}$-points of any simple Lie algebra of classical type. We then get the following:

THEOREM 1.7.   *Let $L$ be a ring additively isomorphic to $\mathbb{Z}^d$. Then there exist some rational number $\alpha(L) \in \mathbb{Q}$, a nonnegative integer $b(L) \in \mathbb{N}$ and some real numbers $c, c' \in \mathbb{R}$ such that $\zeta_L(s)$ has abscissa of convergence $\alpha(L)$ and*

$$
\begin{aligned}
s_N(L) &:= a_1(L) + a_2(L) + \cdots + a_N(L) \sim c \cdot N^{\alpha(L)} \left(\log N\right)^{b(L)}, \\
s_N^{\alpha(L)}(L) &:= a_1(L) + a_2(L) 2^{-\alpha(L)} + \cdots + a_N(L) N^{-\alpha(L)} \sim c \cdot (\log N)^{b(L)+1}
\end{aligned}
$$

*for $N \to \infty$.*

There is a similar theorem for the invariant $a_n^{\triangleleft}(L)$ counting ideals.

We actually prove this theorem as part of our proof of Theorem 1.1, making use of the fact that for a nilpotent group $G$ there is a Lie algebra $L(G)$ defined over $\mathbb{Z}$ with the property that for almost all primes $p$

$$\zeta_{G,p}(s) = \zeta_{L(G),p}(s).$$

This fact was established in [24]. We also use the fact that for those finite number of primes for which this identity does not hold, we still know that $\zeta_{G,p}(s)$ is a rational function whose abscissa of convergence coincides with that of $\zeta_{L(G),p}(s)$.

In [23] the first author and Ph.D. student Gareth Taylor have calculated the zeta function of the Lie algebra $\mathfrak{sl}_2(\mathbb{Z})$ by performing three blow-ups on the associated polynomial $F_{\mathfrak{sl}_2(\mathbb{Z})}$. The paper shows that our method can even be applied to bad primes ($p = 2$ for $\mathfrak{sl}_2(\mathbb{Z})$) where the resolution of singularities of $F_{\mathfrak{sl}_2(\mathbb{Z})}$ does not have good reduction. It is established in [23] that

$$\zeta_{\mathfrak{sl}_2(\mathbb{Z})}(s) = \zeta(s)\zeta(s-1)\zeta(2s-2)\zeta(2s-1)\zeta(3s-1)^{-1} \cdot \frac{(1 + 6 \cdot 2^{-2s} - 8 \cdot 2^{-3s})}{(1 - 2 \cdot 2^{-3s})}.$$

Note that this example has a single pole at $s = 2$. This means then that the subalgebra growth, in contrast to the 3-dimensional Heisenberg-Lie algebra, is $s_N(\mathfrak{sl}_2(\mathbb{Z})) \sim c \cdot N^2$ for $N \to \infty$ where $c = \frac{20}{31} \cdot \frac{\zeta(2)^2\zeta(3)}{\zeta(5)}$. (This example for good primes had been calculated previously in [10] using work of Ishai Ilani [27]. However the calculations of Ilani are heavy. The simplicity of the calculation in [23] is a good advertisement for the practical value of the methods developed in the current paper.)

(3) Let $G$ be a linear algebraic group over $\mathbb{Q}$. Let $\rho : G \to \mathrm{GL}_n$ be a $\mathbb{Q}$-rational representation. Define the '*local zeta function of the algebraic group $G$ at the representation $\rho$ and the prime $p$*' to be

$$Z_{G,\rho,p}(s) = \int_{G^+} |\det \rho(g)|^s \, \mu_G(g)$$

where $G^+ = \rho^{-1}\left(\rho\left(G\left(\mathbb{Q}_p\right)\right) \cap \mathrm{M}_n\left(\mathbb{Z}_p\right)\right)$ and $\mu_G$ denotes the right Haar measure on $G(\mathbb{Q}_p)$ normalized such that $\mu_G\left(G\left(\mathbb{Z}_p\right)\right) = 1$.

We define the '*global zeta function of $G$ at the representation $\rho$*' to be the Euler product

$$Z_{G,\rho}(s) = \prod_p Z_{G,\rho,p}(s).$$

Such zeta functions were first studied by Hey and Tamagawa in the case that $G = \mathrm{GL}_{l+1}$ where $Z_{G,\rho}(s) = \zeta(s) \cdots \zeta(s-l)$. Note that this zeta function is precisely the zeta function $\zeta_{\mathbb{Z}^{l+1}}(s)$. More generally in [24] the zeta functions $Z_{G,\rho}(s)$ are shown to count subgroups $H$ in a nilpotent group $\Gamma$ with the property that the profinite completions are isomorphic; i.e. $\widehat{H} \cong \widehat{\Gamma}$. In this case the algebraic group $G$ is the automorphism group of $\Gamma$. A result of Bryant and Groves shows that any algebraic group can be realised modulo a unipotent group as the automorphism group of a nilpotent group. In [15] an explicit expression is given for the local factors of a class of nilpotent groups in terms of the combinatorics of the building of the algebraic group. The local zeta functions $Z_{G,\rho,p}(s)$ can be expressed in terms of cone integrals. Hence our results apply to these zeta functions.

Although our results imply we can meromorphically continue the zeta function $Z_{G,\rho}(s)$ past its abscissa of convergence, this zeta function in general has a natural boundary, except for the case of $G = \mathrm{GL}_{l+1}$ (see [13]). However we have discovered a procedure which produces something we call the *ghost zeta function* associated to $Z_{G,\rho}(s)$ which often turns out to have a meromorphic continuation to the whole complex plane (see [16] and [17]).

(4) Let $g(n, c, d)$ be the number of finite nilpotent groups of size $n$ of class bounded by $c$ and generated by at most $d$ elements. In [14] the zeta function $\zeta_{\mathcal{N}(c,d)} = \sum_{n=1}^{\infty} g(n, c, d)n^{-s}$ is shown to be expressible as the Euler product of $p$-adic cone integrals. Hence the results of this paper imply that asymptotically $g(n, c, d)$ behaves as follows:

$$g(1, c, d) + g(2, c, d) + \cdots + g(N, c, d) \sim c \cdot N^{\alpha} \left(\log N\right)^b$$

for $N \to \infty$ where $\alpha \in \mathbb{Q}$, $b \in \mathbb{N}$ and $c \in \mathbb{R}$. The details are explained in [9] and [14].

(5) The Igusa zeta function of a polynomial $f(\mathbf{x})$ is defined as

$$Z(s) = \int_{\mathbb{Z}_p^m} |f(\mathbf{x})|^s \, |d\mathbf{x}| \, .$$

Hence it is a particular example of a cone integral where the cone condition $\psi$ is empty. The global zeta function that one can define as the Euler product of these Igusa zeta functions (normalized to have constant coefficient 1) is a special case of our analysis. We consider in a future paper [20] the analytic properties of such global Igusa zeta functions and in particular that they appear to have natural boundaries in a similar fashion to the examples discussed in [13]. In [34] Ono considered a special case of these global Igusa zeta functions and established their region of convergence. He considers the case where the polynomial $f(\mathbf{x})$ is absolutely irreducible and makes use of the Lang-Weil inequality on the number of rational points of a variety as we have. In the special case that the hyper-surface $f(\mathbf{x}) = 0$ is nonsingular, he demonstrates some analytic continuation. Our work may be seen as a vast generalization of Ono's results.

The results of this paper were previously announced in [18].

## Notation

$\mathbb{Q}_p$ denotes the field of $p$-adic numbers.

$\mathbb{Z}_p$ denotes the ring of $p$-adic integers.

For $x \in \mathbb{Q}_p$, $|x|$ denotes $p^{-v(x)}$ where $v(x)$ is the $p$-adic valuation of $x$.

$\mathbb{N}$ denotes the set $\{0, 1, 2, \dots\}$.

$\mathbb{N}_{>0}$ denotes the set $\{1, 2, \dots\}$.

$\mathbb{R}_{>0}$ denotes the set $\{s \in \mathbb{R} : s > 0\}$.

$\mathbb{R}_{\geq 0}$ denotes the set $\{s \in \mathbb{R} : s \geq 0\}$.

$\mathbb{Z}_p^*$ denotes the units of $\mathbb{Z}_p$.

$f(n) \sim g(n)$ means $f(n)/g(n)$ tends to 1 as $n$ tends to infinity.

*Acknowledgements.* We should like to thank Jürgen Elstrodt for discussions concerning the Tauberian theorem. We also thank Benjamin Klopsch and Dan Segal for alerting us to the potential dangers of bad primes in applying the Tauberian theorem. The first author would like to thank the Royal Society, the Max-Planck-Institute in Bonn and the Heinrich Heine Universität in Düsseldorf for support and hospitality during the preparation of this paper.

## 2. An explicit formula for cone integrals

In this section we give a proof of Theorem 1.3 and recall from the introduction the definition of a cone integral:

*Definition* 2.1. (1) Call a formula $\psi(\mathbf{x})$ in the first order language for the valued field $\mathbb{Q}_p$ *a cone condition over* $\mathbb{Q}$ if there exist nonzero polynomials $f_i(\mathbf{x}), g_i(\mathbf{x})(i = 1, \dots, l)$ over $\mathbb{Q}$ in the variables $\mathbf{x} = x_1, \dots, x_m$ such that $\psi(\mathbf{x})$ is a conjunction of formulas

$$v(f_i(\mathbf{x})) \leq v(g_i(\mathbf{x}))$$

for $i = 1, \dots, l$.

(2) Given a cone condition $\psi(\mathbf{x})$ over $\mathbb{Q}$ and nonzero polynomials $f_0$ and $g_0$ with coefficients in $\mathbb{Q}$, we call an integral

$$Z_{\mathcal{D}}(s, p) = \int_{V_p = \left\{\mathbf{x} \in \mathbb{Z}_p^m : \psi(\mathbf{x}) \text{ is valid}\right\}} |f_0(\mathbf{x})|^s \, |g_0(\mathbf{x})| \, |dx|$$

*a cone integral defined over* $\mathbb{Q}$, where $|dx|$ is the normalized additive Haar measure on $\mathbb{Z}_p^m$ and $\mathcal{D} = \{f_0, g_0, f_1, g_1, \dots, f_l, g_l\}$ is called the *cone integral data*.

We are going to use resolution of singularities to get an explicit formula for such cone integrals valid for almost all primes $p$. We follow Section 5 of [5].

*Definition* 2.2. A resolution $(Y, h)$ for a polynomial $F$ over $\mathbb{Q}$ consists of a closed integral subscheme $Y$ of $\mathbf{P}^k_{X_{\mathbb{Q}}}$ (where $X_{\mathbb{Q}} = \mathrm{Spec}(\mathbb{Q}[\mathbf{x}])$ and $\mathbf{P}^k_{X_{\mathbb{Q}}}$ denotes projective $k$-space over the scheme $X_{\mathbb{Q}}$) and the morphism $h : Y \to X$ which is the restriction to $Y$ of the projection morphism $\mathbf{P}^k_{X_{\mathbb{Q}}} \to X_{\mathbb{Q}}$, such that

(i) $Y$ is smooth over $\mathrm{Spec}(\mathbb{Q})$;

(ii) the restriction $h : Y \backslash h^{-1}(D) \to X \backslash D$ is an isomorphism (where $D = \mathrm{Spec}\left(\frac{\mathbb{Q}[\mathbf{x}]}{(F)}\right) \subset X_{\mathbb{Q}}$); and

(iii) the reduced scheme $(h^{-1}(D))_{\mathrm{red}}$ associated to $h^{-1}(D)$ has only normal crossings (as a subscheme of $Y$).

Let $E_i$, $i \in T$, be the irreducible components of the reduced scheme $(h^{-1}(D))_{\mathrm{red}}$ over $\mathrm{Spec}(\mathbb{Q})$. For $i \in T$, let $N_i$ be the multiplicity of $E_i$ in the divisor of $F \circ h$ on $Y$ and let $\nu_i - 1$ be the multiplicity of $E_i$ in the divisor of $h^*(dx_1 \wedge \cdots \wedge dx_m)$. The $(N_i, \nu_i)$ $i \in T$, are called the *numerical data* of the resolution $(Y, h)$ for $F$.

Let us recall some necessary facts about reduction of varieties mod $p$. When $X = X_{\mathbb{Q}} = \mathrm{Spec}(\mathbb{Q}[\mathbf{x}])$ one defines the reduction mod $p$ of a closed integral subscheme $Y$ of $\mathbf{P}^k_{X_{\mathbb{Q}}}$ as follows: let $\widetilde{X} = \mathrm{Spec}(\mathbb{Z}[\mathbf{x}])$ and $\widetilde{Y}$ be the scheme-theoretic closure of $Y$ in $\mathbf{P}^k_{\widetilde{X}}$. Then the reduction mod $p$ of $Y$ is the scheme $\widetilde{Y} \times_{\mathbb{Z}} \mathrm{Spec}(\mathbb{F}_p)$ and we denote it by $\overline{Y}$. Let $\widetilde{h} : \widetilde{Y} \to \widetilde{X}$ be the restriction to $\widetilde{Y}$ of the projection morphism $\mathbf{P}^k_{\widetilde{X}} \to \widetilde{X}_{\mathbb{Q}}$ and $\overline{h} : \overline{Y} \to \overline{X}$ be obtained from $\widetilde{h}$ by base extension.

*Definition* 2.3. A resolution $(Y, h)$ for $F$ over $\mathbb{Q}$ has *good reduction* mod $p$ if

(1) $\overline{Y}$ is smooth over $\mathrm{Spec}(\mathbb{F}_p)$;

(2) $\overline{E_i}$ is smooth over $\mathrm{Spec}(\mathbb{F}_p)$, for each $i \in T$, and $\bigcup_{i \in T} \overline{E_i}$ has only normal crossings as a subscheme of $\overline{Y}$; and

(3) $\overline{E_i}$ and $\overline{E_j}$ have no common irreducible components, when $i \neq j$.

Note that a resolution over $\mathbb{Q}$ has good reduction for almost all primes $p$ (see Theorem 2.4 of [5]).

Let $(Y^o, h^o)$ be a resolution for the polynomial $F = \prod_{i=0}^l f_i \cdot g_i$ over $\mathbb{Q}$, and $p$ be any prime such that $(Y^o, h^o)$ has good reduction mod $p\mathbb{Z}_p$ and $\overline{\prod_{i=0}^l f_i \cdot g_i} \neq 0$. Here $\overline{\phantom{-}}$ means reduction mod $p$ . Let $(Y, h)$ be the resolution over $\mathbb{Q}_p$ obtained from $(Y^o, h^o)$ by base extension.

Let $a \in \overline{Y}(\mathbb{F}_p)$. Since we consider $\overline{Y}$ as a closed subscheme of $\widetilde{Y}$, $a$ is also a closed point of $\widetilde{Y}$. Let $T_a = \left\{ i \in T : a \in \overline{E_i} \right\} = \left\{ i \in T : a \in \widetilde{E_i} \right\}$. Let $r = \operatorname{card} T_a$ and $T_a = \{i_1, \ldots, i_r\}$. Then in the local ring $O_{\widetilde{Y},a}$ we can write

$$F \circ \widetilde{h} = u c_1^{N_{i_1}} \cdots c_r^{N_{i_r}}$$

where $c_i \in O_{\widetilde{Y},a}$ generates the ideal of $\widetilde{E_{i_j}}$ in $O_{\widetilde{Y},a}$ and $u$ is a unit in $O_{\widetilde{Y},a}$. Since $f_i$ and $g_i$ divide $F$ we can also write for $i = 0, \ldots, l$

$$
\begin{aligned}
f_i \circ \widetilde{h} &= u_{f_i} c_1^{N_{i_1}(f_i)} \cdots c_r^{N_{i_r}(f_i)}, \\
g_i \circ \widetilde{h} &= u_{g_i} c_1^{N_{i_1}(g_i)} \cdots c_r^{N_{i_r}(g_i)}.
\end{aligned}
$$

Put

$$J_a(s,p) = \int_{\theta^{-1}(a) \cap h^{-1}(V_p)} |f_0 \circ h|^s \, |g_0 \circ h| \, |h^*(dx_1 \wedge \cdots \wedge dx_m)|$$

where we define $\theta$ as follows: Let $H = \{b \in Y(\mathbb{Q}_p) : h(b) \in \mathbb{Z}_p^m\}$. A point $b \in H \subset Y(\mathbb{Q}_p)$ can be represented by its coordinates $(x_1, \ldots, x_m, y_0, \ldots, y_k)$ in $\mathbb{Q}_p^m \times \mathbf{P}_X^k(\mathbb{Q}_p)$ where $(x_1, \ldots, x_m) \in \mathbb{Z}_p^m$ and $y_0, \ldots, y_k$ are homogeneous coordinates which can therefore be chosen such that $\min_{i=0,\ldots,k} \operatorname{ord} y_i = 0$. The map $\theta : H \to \overline{Y}(\mathbb{F}_p)$ is then defined as follows: $\theta(b) = (\overline{x_1}, \ldots, \overline{x_m}, \overline{y_0}, \ldots, \overline{y_k}) \in \overline{Y}(\mathbb{F}_p) \subset \mathbf{P}_{\overline{X}}^k(\mathbb{F}_p)$. Then $Z_{\mathcal{D}}(s,p) = \sum_{a \in \overline{Y}(\mathbb{F}_p)} J_a(s,p)$. Now we have

$$J_a(s,p) = \int_{\theta^{-1}(a) \cap h^{-1}(V_p)} |c_1|^{N_{i_1}(f_0)s + N_{i_1}(g_0) + \nu_{i_1} - 1} \cdots$$
$$\cdots |c_r|^{N_{i_r}(f_0)s + N_{i_r}(g_0) + \nu_{i_r} - 1} |dc_1 \wedge \cdots \wedge dc_m|.$$

Since $\overline{c_1}, \ldots, \overline{c_m}$ belong to the maximal ideal of $O_{\overline{Y},a}$, we have $c_1(b), \ldots, c_m(b) \in p\mathbb{Z}_p$ for all $b \in \theta^{-1}(a)$. The map

$$
\begin{aligned}
c &: & \theta^{-1}(a) &\to (p\mathbb{Z}_p)^m \\
& & b &\mapsto (c_1(b), \ldots, c_m(b))
\end{aligned}
$$

is a bijection. Hence

(2.1)
$$J_a(s,p) = \int_{V_p'} |y_1|^{N_{i_1}(f_0)s + N_{i_1}(g_0) + \nu_{i_1} - 1} \cdots |y_r|^{N_{i_r}(f_0)s + N_{i_r}(g_0) + \nu_{i_r} - 1} |dy_1| \cdots |dy_m|$$

where $V_p'$ is the set of all $y = (y_1, \ldots, y_m) \in (p\mathbb{Z}_p)^m$ satisfying, for $i = 1, \ldots, l$,

$$\sum_{j=1}^{r} N_{i_j}(f_i) \operatorname{ord}(y_j) \leq \sum_{j=1}^{r} N_{i_j}(g_i) \operatorname{ord}(y_j).$$

Let $A_{j,a} = N_{i_j}(f_0)$ and $B_{j,a} = N_{i_j}(g_0) + \nu_{i_j}$ for $j = 1, \ldots, r$ and $A_{j,a} = 0, B_{j,a} = 1$ for $j > r$. Then

$$
\begin{aligned}
J_a(s,p) &= \sum_{(k_1,\ldots,k_m)\in\Lambda} p^{-\sum_{j=1}^m k_j(A_{j,a}s+B_{j,a}-1)}(p^{-k_1}-p^{-k_1-1})\cdots \\
&\qquad\qquad\qquad\qquad\qquad\qquad \cdots(p^{-k_m}-p^{-k_m-1}) \\
&= (1-p^{-1})^m \sum_{(k_1,\ldots,k_m)\in\Lambda} p^{-\sum_{j=1}^m k_j(A_{j,a}s+B_{j,a})}
\end{aligned}
$$

where $A_{j,a} \in \mathbb{N}$ and $B_{j,a} \in \mathbb{N}$ and

$$
\Lambda = \left\{ (k_1,\ldots,k_m) \in \mathbb{N}_{>0}^m : \sum_{j=1}^r N_{i_j}(f_i)k_j \leq \sum_{j=1}^r N_{i_j}(g_i)k_j \text{ for } i = 1,\ldots,l \right\}.
$$

Thus $\Lambda$ is the intersection of $\mathbb{N}_{>0}^m$ and a rational convex polyhedral cone $C$ in $\mathbb{R}_{>0}^m$. We can write this cone as a disjoint union of simplicial cones $C_1,\ldots,C_w$ of the form:

$$
C_i = \{\alpha_1 v_{i1} + \cdots + \alpha_{m_i} v_{im_i} : \alpha_j \in \mathbb{R}_{>0}, \text{ for } j = 1,\ldots,m_i\}
$$

where $\{v_{i1},\ldots,v_{im_i}\}$ is a linearly independent set of vectors in $\mathbb{R}^m$ with non-negative integer coordinates and with the property that a fundamental region of the lattice spanned by $v_{i1},\ldots,v_{im_i}$ has no lattice point of $\mathbb{Z}^m$ in its interior (see p. 123–124 of [1]). Then $\Lambda$ can be written as the disjoint union of the following sets:

$$
\Lambda_i = \{l_1 v_{i1} + \cdots + l_{m_i} v_{im_i} : l_j \in \mathbb{N}_{>0} \text{ for } j = 1,\ldots,m_i\}.
$$

Put $v_{jk} = (q_{jk1},\ldots,q_{jkm}) \in \mathbb{N}^m$ for $k = 1,\ldots,m_j$. Hence

$$
J_a(s,p) = (1-p^{-1})^m \sum_{j=1}^w \prod_{k=1}^{m_j} \frac{p^{-(A_{k,a,j}s+B_{k,a,j})}}{1-p^{-(A_{k,a,j}s+B_{k,a,j})}}
$$

where $A_{k,a,j} = \sum_{i=1}^m q_{jki} A_{i,a} \in \mathbb{N}$ and $B_{k,a,j} = \sum_{i=1}^m q_{jki} B_{i,a} \in \mathbb{N}$.

Notice that the above calculations just depended on which components $\overline{E_i}$ contained $a$. If $T_{a_1} = T_{a_2}$ then $J_{a_1}(s,p) = J_{a_2}(s,p)$. So for each $I \subset T$ let

$$
c_{p,I} = \operatorname{card}\{a \in \overline{Y}(\mathbb{F}_p) : a \in \overline{E_i} \text{ if and only if } i \in I\}
$$

and put $A_{k,I,j} = A_{k,a,j}$ and $B_{k,I,j} = B_{k,a,j}$ for any $a \in \{a \in \overline{Y}(\mathbb{F}_p) : a \in \overline{E_i}$ if and only if $i \in I\}$ where $j = 1,\ldots,w_I$ and $w_I$ is the number of simplicial cones defined by the linear inequalities corresponding to $I$. Then we have a final formula for $Z_{\mathcal{D}}(s,p)$:

$$
(2.2) \qquad Z_{\mathcal{D}}(s,p) = (1-p^{-1})^m \sum_{I \subset T} c_{p,I} \sum_{j=1}^{w_I} \prod_{k=1}^{m_j} \frac{p^{-(A_{k,I,j}s+B_{k,I,j})}}{1-p^{-(A_{k,I,j}s+B_{k,I,j})}}.
$$

Note that if $A_{k,I,j} = 0$ and $B_{k,I,j} = 1$, which will correspond to a bit of the integral like $\int_{p\mathbb{Z}_p} |dy_m|$, then we get $(1 - p^{-1}) \cdot \frac{p^{-1}}{(1-p^{-1})}$ which is correct.

This completes the proof of Theorem 1.3. Note that this expression (2.2) for $Z_{\mathcal{D}}(s, p)$ holds for all primes for which the resolution $(Y, h)$ had good reduction.

We could also consider a cone integral defined over $\mathbb{Q}_p$ rather than $\mathbb{Q}$ whose cone data $\mathcal{D}$ consisted of polynomials in $\mathbb{Q}_p[\mathbf{x}]$. Our formula (2.2) would still be valid for such integrals provided that the resolution had good reduction mod $p$.

Notice that, as we vary $p$, the only things in this formula which depend on $p$ are the terms $c_{p,I}$.

We should note that there is one term which is always a constant term in the expression for our final formula (2.2) corresponding to the subset $I = \emptyset$; then $w_\emptyset = 1$, $m_1 = m$ and $A_{k,\emptyset,1} = 0$ and $B_{k,\emptyset,1} = 1$ for $k = 1, \ldots, m$. Hence the term corresponding to the subset $I = \emptyset$ has the following form:

$$(2.3)$$
$$(1 - p^{-1})^m c_{p,\emptyset} \sum_{j=1}^{w_\emptyset} \prod_{k=1}^{m_j} \frac{p^{-(A_{k,\emptyset,j}s + B_{k,\emptyset,j})}}{1 - p^{-(A_{k,\emptyset,j}s + B_{k,\emptyset,j})}} = c_{p,\emptyset}(1 - p^{-1})^m \frac{p^{-m}}{(1 - p^{-1})^m}$$
$$= c_{p,\emptyset} p^{-m}.$$

Since the restriction $h : Y \backslash h^{-1}(D) \to X \backslash D$ is an isomorphism (where $D = \mathrm{Spec}\left(\frac{\mathbb{Q}[\mathbf{x}]}{(F)}\right) \subset X_{\mathbb{Q}}$)

$$c_{p,\emptyset} = \mathrm{card}\{a \in \overline{Y}(\mathbb{F}_p) : a \notin \overline{E_i} \text{ for all } i \in T\}$$
$$= \mathrm{card}\,\overline{X}(\mathbb{F}_p) - \mathrm{card}\,\overline{D}(\mathbb{F}_p).$$

The term (2.3) is part of the constant term of the rational function $Z_{\mathcal{D}}(s, p)$. The other parts of the constant term come from those $I \subset T$ and $j \in \{1, \ldots, w_I\}$ such that $A_{k,I,j} = 0$ for all $k = 1, \ldots, m_j$.

Note that by dimension arguments for $p$ large enough, $\mathrm{card}\,\overline{X}(\mathbb{F}_p) > \mathrm{card}\,\overline{D}(\mathbb{F}_p)$. Hence $c_{p,0} > 0$ for almost all primes $p$ and the constant term $a_{p,0}$ in a cone integral is nonzero for almost all primes $p$ as promised in the introduction. We give a lower bound for this constant in Section 4.

## 3. A second explicit expression for cone integrals

The explicit expression (2.2) determined in the previous section has a number of advantages. It expresses the function as a sum over the subsets of $I$ which identifies precisely the bits $c_{p,I}$ which depend on $p$. This form of the sum is also more amenable to Denef and Meuser's proof that the Igusa local zeta function (where $\psi$ is the empty condition) satisfies a functional equation (see [6]).

However, for the analysis of the analytic properties of the global zeta function, as explained in the introduction, it is preferable to work with a second explicit formula (to be established) where the cone integrals are written as a sum over open simplicial pieces of a single cone defined in card $T$ dimensions, where each open simplicial piece of the cone gets a weight according to the size of $I$ and $c_{p,I}$.

We give a proof in this section of Theorem 1.4 where all the data in the formula, e.g. $A_j$ and $B_j$, are identified explicitly in terms of the numerical data of the resolution and the underlying cone.

The cone is defined as follows:

$$\overline{D_T} = \left\{ (x_1, \ldots, x_t) \in \mathbb{R}^t_{\geq 0} : \sum_{j=1}^t N_j(f_i)x_j \leq \sum_{j=1}^t N_j(g_i)x_j \text{ for } i = 1, \ldots, l \right\}$$

where card $T = t$ and $\mathbb{R}_{\geq 0} = \{x \in \mathbb{R} : x \geq 0\}$; so this is a closed cone. Denote the lattice points in $\overline{D_T}$ by $\overline{\Delta_T}$, i.e. $\overline{\Delta_T} = \overline{D_T} \cap \mathbb{N}^t$. We can write this cone as a disjoint union of open simplicial pieces called $R_k$, $k = 0, 1, \ldots, w$ where a fundamental region for the lattice points of $R_k$ has no lattice points in its interior. We shall assume that $R_0 = (0, \ldots, 0)$ and that the next $q$ pieces are all the open one-dimensional edges in our choice of simplicial decomposition for the cone $\overline{D_T}$: for $k = 1, \ldots, q$,

$$R_k = \{\alpha \mathbf{e}_k = \alpha(q_{k1}, \ldots, q_{kt}) : \alpha > 0\}.$$

Since these are all the one-dimensional edges, for any $k \in \{0, \ldots, w\}$ there exists some subset $M_k \subset \{1, \ldots, q\}$ such that

$$R_k = \left\{ \sum_{j \in M_k} \alpha_j \mathbf{e}_j : \alpha_j > 0 \text{ for all } j \in M_k \right\}.$$

Note that $m_k := \text{card } M_k \leq t$.

Define for each $k = 1, \ldots, q$ the following constants:

$$(3.1) \qquad \begin{aligned} A_k &= \sum_{j=1}^t q_{kj} N_j(f_0), \\ B_k &= \sum_{j=1}^t q_{kj} \left( N_j(g_0) + \nu_j \right). \end{aligned}$$

For each subset $I \subset T$ we previously defined a rational convex polyhedral cone $C_I$ with lattice points $\Lambda_I$ which we broke down into simplicial cones $C_1^I, \ldots, C_{w_I}^I$ with corresponding lattice points $\Lambda_1^I, \ldots, \Lambda_{w_I}^I$. These were cones in the open positive quadrant $\mathbb{R}^m_{>0}$. We are going to use the new cone $\overline{D_T}$ to express the same rational function that we associated to $C_I$.

For each $I \subset T$ define:

$$D_I = \left\{ (k_1, \ldots, k_t) \in \overline{D_T} : k_i > 0 \text{ if } i \in I \text{ and } k_i = 0 \text{ if } i \in T \backslash I \right\},$$
$$\Delta_I = D_I \cap \mathbb{N}^t.$$

So $\overline{D_T} = \bigcup_{I \subset T} D_I$, a disjoint union. The reason we chose the notation $\overline{D_T}$ is of course because it is the closure then of $D_T$ defined as above.

For each $I \subset T$ there is then a subset $W_I \subset \{0, \ldots, w\}$ so that

$$D_I = \bigcup_{k \in W_I} R_k.$$

We now have the following:

LEMMA 3.1.

$$(1 - p^{-1})^m c_{p,I} \sum_{j=1}^{w_I} \prod_{k=1}^{m_j} \frac{p^{-(A_{k,I,j}s + B_{k,I,j})}}{1 - p^{-(A_{k,I,j}s + B_{k,I,j})}}$$
$$= c_{p,I} \sum_{k \in W_I} (1 - p^{-1})^{|I|} p^{-(m-|I|)} \prod_{j \in M_k} \frac{p^{-(A_j s + B_j)}}{1 - p^{-(A_j s + B_j)}}.$$

COROLLARY 3.2. *Set* $c_{p,k} = c_{p,I}$ *and* $I_k = I$ *if* $k \in W_I$. *Then for all primes* $p$ *for which the resolution has good reduction,*

$$Z_{\mathcal{D}}(s, p) = \sum_{k=0}^{w} (p-1)^{|I_k|} p^{-m} c_{p,k} \prod_{j \in M_k} \frac{p^{-(A_j s + B_j)}}{1 - p^{-(A_j s + B_j)}}.$$

*Proof.* We go back to the integral expression (2.1) for $J_a(s, p)$ where $a \in \{a \in \overline{Y}(\mathbb{F}_p) : a \in \overline{E_i} \text{ if and only if } i \in I\}$. The calculations of the previous section gave:

$$J_a(s, p) = (1 - p^{-1})^m \sum_{j=1}^{w_I} \prod_{k=1}^{m_j} \frac{p^{-(A_{k,I,j}s + B_{k,I,j})}}{1 - p^{-(A_{k,I,j}s + B_{k,I,j})}}.$$

We can rewrite the expression for $J_a(s, p)$ as

$$J_a(s, p) = p^{-(m-|I|)} \int_{V_p'} \prod_{i \in I} |z_i|^{N_i(f_0)s + N_i(g_0) + \nu_i - 1} \prod_{i \in I} |dz_i|$$

where $V_p'$ is now the set of $(z_i)_{i \in I} \in (p\mathbb{Z}_p)^{|I|}$ satisfying for $j = 1, \ldots, l$

$$\sum_{i \in I} N_i(f_j) \mathrm{ord}(z_i) \leq \sum_{i \in I} N_i(g_j) \mathrm{ord}(z_i).$$

By the definition of $\Delta_I = D_I \cap \mathbb{N}^t \subset \overline{D_T}$ this then reduces to

$$
\begin{aligned}
J_a(s,p) &= p^{-(m-|I|)}(1-p^{-1})^{|I|} \sum_{(k_1,\ldots,k_t) \in \Delta_I} p^{-\sum_{j=1}^{t} k_j(N_j(f_0)s + N_j(g_0) + \nu_j)} \\
&= \sum_{k \in W_I} p^{-(m-|I|)}(1-p^{-1})^{|I|} \\
&\quad \times \sum_{(k_1,\ldots,k_t) \in R_k \cap \mathbb{N}^t} p^{-\sum_{j=1}^{t} k_j(N_j(f_0)s + N_j(g_0) + \nu_j)}
\end{aligned}
$$

since $D_I = \bigcup_{k \in W_I} R_k$. Because

$$
R_k \cap \mathbb{N}^t = \left\{ \sum_{j \in M_k} \alpha_j \mathbf{e}_j : \alpha_j \in \mathbb{N}_{>0} \text{ for all } j \in M_k \right\},
$$

by making a change of variable as in the previous section and using the definitions (3.1) of $A_k$ and $B_k$ for $k = 1, \ldots, q$, we get

$$
J_a(s,p) = \sum_{k \in W_I} p^{-(m-|I|)}(1-p^{-1})^{|I|} \prod_{j \in M_k} \frac{p^{-(A_j s + B_j)}}{1 - p^{-(A_j s + B_j)}}.
$$

This completes the proof of the lemma. $\qquad\qquad\square$

Recall that even if the dimension of the simplicial piece has gone down (i.e. $m_k < |I|$ ), we will still get a $(1-p^{-1})$ for each variable. For example, the integral $\int_{\mathbb{Z}_p^2} |x|^s |y|^s$ over $v(x) = v(y)$ is $(1-p^{-1})^2(1-p^{-2(s+1)})^{-1}$. The second point to note is that $c_{p,I} = 0$ for any $I \subset T$ for which card $I > m$, where $m$ is the number of variables in the original integral.

We conclude this section by showing that even for primes of bad reduction, the rational expression for these local factors is not too far from the explicit expressions established here for good primes. In particular we can establish that the candidate local poles for the bad primes are a subset of the candidate poles $\{-B_j/A_j : j = 1, \ldots, q\}$ for the expressions for the good primes. We follow Igusa's original proof of the rationality of the local zeta functions (see [25] or the more recent volume [26]). The essential observation is that the resolution of singularities over $\mathbb{Q}$ is still a resolution of singularities for all $\mathbb{Q}_p$ regardless of whether the prime has good or bad reduction. We can then establish the following:

PROPOSITION 3.3. *For each prime $p$, there exists a finite set $\mathcal{B}_p$ such that for each $b \in \mathcal{B}_p$ there is an associated subset $I_b \subset T$ and integer $e_b$ such that*

$$
(3.2) \quad Z_{\mathcal{D}}(s,p) = \sum_{b \in \mathcal{B}_p} \left( \sum_{k \in W_{I_b}} p^{-(m-|I_b|)}(1-p^{-1})^{|I_b|} \prod_{j \in M_k} \frac{p^{-e_b(A_j s + B_j)}}{1 - p^{-(A_j s + B_j)}} \right).
$$

*Proof.* Since $(Y, h)$ is a resolution of singularities over $\mathbb{Q}_p$ for the polynomial $F$, $H = \{x \in Y(\mathbb{Q}_p) : h(b) \in \mathbb{Z}_p^m\}$ can be written as a finite disjoint union of open subsets $U(b)$, $b \in \mathcal{B}_p$, such that $U(b)$ has a surjective chart $\phi_{U(b)} : U(b) \to p^{e_b}\mathbb{Z}_p^m$ for some $e_b \in \mathbb{N}$. Associated to each chart there is a subset $I_b = \{i_1, \ldots, i_r\} \subset T$ such that for every $y \in U(b)$ with $\phi_{U(b)}(y) = (y_1, \ldots, y_m)$

$$\left| f_i \circ h \circ \phi_{U(b)}^{-1}(y_1, \ldots, y_m) \right| = \left| y_1^{N_{i_1}(f_i)} \right| \cdots \left| y_r^{N_{i_r}(f_i)} \right|,$$

$$\left| g_i \circ h \circ \phi_{U(b)}^{-1}(y_1, \ldots, y_m) \right| = \left| y_1^{N_{i_1}(g_i)} \right| \cdots \left| y_r^{N_{i_r}(g_i)} \right|.$$

We can use these charts as before to express our cone integral for $p$ as

$$Z_{\mathcal{D}}(s, p) = \sum_{b \in \mathcal{B}_p} J_b(s, p)$$

where

$$J_b(s, p) = \int_{V_p'} |y_1|^{N_{i_1}(f_0)s + N_{i_1}(g_0) + \nu_{i_1} - 1} \cdots |y_r|^{N_{i_r}(f_0)s + N_{i_r}(g_0) + \nu_{i_r} - 1} |dy_1| \cdots |dy_m|$$

and $V_p'$ is the set of all $(y_1, \ldots, y_m) \in p^{e_b}\mathbb{Z}_p^m$ satisfying, for $i = 1, \ldots, l$,

$$\sum_{j=1}^{r} N_{i_j}(f_i)\mathrm{ord}(y_j) \leq \sum_{j=1}^{r} N_{i_j}(g_i)\mathrm{ord}(y_j).$$

Our analysis of these integrals, as before, implies the statement of the theorem. □

COROLLARY 3.4. *For all primes $p$, the abscissa of convergence of $Z_{\mathcal{D}}(s, p)$ is one of the rational number $-B_j/A_j$ where $j = 1, \ldots, q$ and $A_j \neq 0$.*

## 4. Zeta functions defined as Euler products of cone integrals

We now turn to analysing the global behaviour of a product of these cone integrals over all primes $p$.

We make some normalisation of the cone integrals so that the constant coefficient of the local factors is 1.

*Definition* 4.1. Let $Z_{\mathcal{D}}(s, p)$ be a cone integral defined over $\mathbb{Q}$. Then denote by $a_{p,0}$ *the constant coefficient of the power series in $p^{-s}$ representing the rational function $Z_{\mathcal{D}}(s, p)$.*

*Definition* 4.2. We say that a function $Z(s)$ is *an Euler product of cone integrals over $\mathbb{Q}$ with cone integral data $\mathcal{D}$* if

$$Z(s) = Z_{\mathcal{D}}(s) = \prod_{p \text{ prime, } a_{p,0} \neq 0} \left( a_{p,0}^{-1} \cdot Z_{\mathcal{D}}(s, p) \right).$$

*Remark* 1. In our analysis we shall want to exclude the trivial case where $Z_{\mathcal{D}}(s,p) = a_{p,0}$ for all $p$. A constant function converges everywhere so that its analysis is not of interest. To check whether a cone integral is constant we need to check whether

$$\left\{ \mathbf{x} \in \mathbb{Z}_p^m : f_0(\mathbf{x}) = 0 \bmod p \text{ and } \psi(\mathbf{x}) \text{ is valid} \right\}$$

is empty for all $p$.

THEOREM 4.3. *A nonconstant function $Z(s)$ defined as an Euler product of cone integrals over $\mathbb{Q}$ has a rational abscissa of convergence $\alpha \in \mathbb{Q}$.*

*Proof.* Let $Q = Q_1 \cup Q_2$ denote the finite set of primes where the $Q_1$ are those primes $p$ with bad reduction for which $a_{p,0} \neq 0$ and the $Q_2$ are those $p$ for which $a_{p,0} = 0$. Let $W_k'$ denote the set of those $k$ for which $\sum_{j \in M_k} A_j \neq 0$. Put $W' = \bigcup_{I \subset T} W_I'$. Then:

$$(4.1)$$

$$Z_{\mathcal{D}}(s) = P(s) \cdot \prod_{p \notin Q} \left( 1 + \sum_{k \in W'} \frac{c_{p,k}}{a_{p,0}} p^{-m} (p-1)^{|I_k|} \prod_{j \in M_k} \frac{p^{-(A_j s + B_j)}}{1 - p^{-(A_j s + B_j)}} \right)$$

where $P(s) = \prod_{p \in Q_1} Z_{\mathcal{D}}(s,p) = \prod_{p \in Q_1} P_p(p^{-s})$ and $P_p(X)$ is a rational function in $\mathbb{Q}(X)$. Note that the abscissa of convergence of $P_p(p^{-s})$ is a rational number since by the formula for bad primes (3.2) the denominator is a product of terms of the form $(1 - p^{-(A_j s + B_j)})$ where $j \in \{1, \ldots, q\}$. It will suffice to prove that

$$\prod_{p \notin Q} \left( 1 + \sum_{k \in W'} \frac{c_{p,k}}{a_{p,0}} p^{-m} (p-1)^{|I_k|} \prod_{j \in M_k} \frac{p^{-(A_j s + B_j)}}{1 - p^{-(A_j s + B_j)}} \right)$$

has a rational abscissa of convergence.

We explain now a few facts about counting points on the reduction of varieties mod $p$.

The Lang-Weil estimate [29] will be a crucial tool:

PROPOSITION 4.4. *There is a constant $C = C(f,k)$ such that every absolutely irreducible quasiprojective variety $E \subset \mathbf{P}^k$ defined over $\mathbb{F}_p$ of degree $f$ and of dimension $d$ satisfies*

$$\left| \operatorname{card} E(\mathbb{F}_p) - p^d \right| \leq C p^{d-1/2}.$$

We shall be interested in counting points on $\overline{X}$, the reduction mod $p$ of a variety $X$ defined over $\mathbb{Q}$ which is irreducible over $\mathbb{Q}$. Let $d$ be the dimension of $X$. The reduction $\overline{X}$ need not be irreducible as a variety over $\mathbb{F}_p$.

Let $\widehat{\mathbb{Q}}$ be the algebraic closure of $\mathbb{Q}$ and $\mathbb{G}$ its Galois group. Consider the decomposition

$$X = X_1 \cup \cdots \cup X_n$$

of $X$ into irreducible components over $\widehat{\mathbb{Q}}$. The Galois group $\mathbb{G}$ acts transitively on the set of components $\{X_1, \ldots, X_n\}$ because $X$ is defined over $\mathbb{Q}$. Since the action is transitive all the $X_i$ also have dimension $d$. Let $U \leq \mathbb{G}$ be the kernel of this action and put $L = \widehat{\mathbb{Q}}^U$. Then $L$ is a finite Galois extension of $\mathbb{Q}$ with Galois group $\mathcal{G} = \mathbb{G}/U$ and every $X_i$ $(i = 1, \ldots, n)$ is defined over $L$.

For every prime number $p$ we choose a prime ideal $\mathfrak{p}$ in $L$ which divides $p$. We write $I_{\mathfrak{p}} \leq D_{\mathfrak{p}} \leq \mathbb{G}$ for the corresponding inertia, respectively decomposition group. If $p$ is unramified in $L$, that is, $I_{\mathfrak{p}} = \{1\}$, we denote by $\mathrm{Frob}_p$ the conjugacy class in $\mathcal{G}$ consisting of the Frobenius elements.

We choose a finite set of primes $S$ such that for every $p \notin S$ the following are satisfied:

(1) the reduction $\overline{X}$ mod $p$ of $X$ is smooth; and

(2) $p$ is unramified in $L$.

We define a function associated to the variety $X$ which will be an essential tool in analysing our zeta functions.

*Definition* 4.5.   Let $X$ be a smooth quasiprojective variety defined over $\mathbb{Q}$ which is irreducible over $\mathbb{Q}$. Define $l_p(X)$ to be the number of irreducible components (defined over $\mathbb{F}_p$) of $\overline{X}$, the reduction mod $p$ of $X$, which are absolutely irreducible. Define

$$V_X(s) = \prod_p \left(1 - l_p(X)p^{-s}\right).$$

We prove the following important properties of this function $V_X(s)$:

LEMMA 4.6.   (1) *The abscissa of convergence of $V_X(s)$ is 1.*

(2) *There is a $\delta > 0$ such that $V_X(s)$ has a meromorphic continuation to the half-plane $\Re(s) > 1 - \delta$.*

*Proof.* We apply the set-up introduced above. The finite group $\mathcal{G}$ acts on $\{X_1, \ldots, X_n\}$ and we write $M$ for the corresponding (complex) permutation module. We obtain a complex finite-dimensional representation $\rho : \mathcal{G} \to \mathrm{GL}(M)$. Denote by $\mathrm{Tr}(\rho(\mathrm{Frob}_p))$ the trace of a representation of the conjugacy class $\mathrm{Frob}_p$ where $p \notin S$. We claim that

$$(4.2) \qquad\qquad \mathrm{Tr}(\rho(\mathrm{Frob}_p)) = l_p(X)$$

for all $p \notin S$. To prove this we let $\overline{X_1}, \ldots, \overline{X_n}$ be the reductions of $X_1, \ldots, X_n$ mod $\mathfrak{p}$. These are absolutely irreducible, smooth quasiprojective varieties defined over the residue field $\mathbb{F}_{\mathfrak{p}}$ corresponding to the prime ideal $\mathfrak{p}$. Let $\mathcal{G}_{\mathfrak{p}}$ be the Galois group of $\mathbb{F}_{\mathfrak{p}}$ over $\mathbb{F}_p$. Then $\mathcal{G}_{\mathfrak{p}}$ acts on the set of components $\overline{X_1}, \ldots, \overline{X_n}$ and the reduction mod $\mathfrak{p}$ is an equivariant map from $\{X_1, \ldots, X_n\}$ to $\left\{\overline{X_1}, \ldots, \overline{X_n}\right\}$ with respect to the natural isomorphism from $D_{\mathfrak{p}}$ to $\mathcal{G}_{\mathfrak{p}}$. Hence the number of fixed points of $\mathcal{G}$ on $\{X_1, \ldots, X_n\}$ is also equal to $l_p(X)$ and formula (4.2) follows.

We are now ready to prove statement (1) of Lemma 4.6. Clearly $l_p(X)$ is bounded by the number $n$ of absolutely irreducible components of $X$. This implies that the abscissa of convergence of $V_X(s)$ is less than or equal to 1. On the other hand the set of primes $p \notin S$ such that $\mathrm{Frob}_p = \{1\}$ has nonzero arithmetic density by Čebotarov's theorem (see [32]). This establishes (1).

Finally to prove (2), we define

$$L_X(s) = \prod_{p \notin S} \det\left(1 - \rho(\mathrm{Frob}_p) \cdot p^{-s}\right)^{-1}.$$

This Euler product is up to finitely many Euler factors equal to the Artin $L$-function of $\rho$. It is well known that $L_X(s)$ has abscissa of convergence 1 and also has a meromorphic continuation to all of $\mathbb{C}$ (see [32]). From formula (4.2) we infer that

$$\det\left(1 - \rho(\mathrm{Frob}_p) \cdot p^{-s}\right) = 1 - l_p(X) \cdot p^{-s} + \sum_{k=2}^{n} a_k p^{-ks}$$

with suitable $a_k$. It follows that

$$\left(1 - l_p(X) \cdot p^{-s}\right)\left(\det\left(1 - \rho(\mathrm{Frob}_p) \cdot p^{-s}\right)^{-1}\right) = 1 + \sum_{k=2}^{\infty} b_k p^{-ks}$$

with $b_k \in \mathbb{C}$. Since the eigenvalues of $\rho(\mathrm{Frob}_p)$ are roots of unity, the $b_k$ can be bounded independently of $p \notin S$. This proves (2). $\qquad\square$

LEMMA 4.7. *Let $X$ be a smooth quasiprojective variety defined over $\mathbb{Q}$ which is irreducible over $\mathbb{Q}$. Let $d$ be the dimension of $X$. There exists $\delta \in \mathbb{R}$ such that for almost all primes $p$*

$$|\mathrm{card}(\overline{X}(\mathbb{F}_p)) - l_p(X)p^d| \leq \delta p^{d-1/2}$$

*and $l_p(X) > 0$ for a dense set of primes.*

*Proof.* Let $Y$ be a smooth quasiprojective irreducible variety over $\mathbb{F}_p$ which is not absolutely irreducible. Then $Y(\mathbb{F}_p)$ is empty. To prove this let $Y = Y_1 \cup \cdots \cup Y_r$ be the decomposition into irreducible components over $\widehat{\mathbb{F}_p}$. The Galois group of $\widehat{\mathbb{F}_p}$ over $\mathbb{F}_p$ acts transitively on the set of components $Y_1, \ldots, Y_r$. So an $\mathbb{F}_p$-point of $Y$ would be contained in every component. If there is more than one component then it is a singular point of $Y$. But $Y$ was assumed to be a smooth variety. Hence no such $\mathbb{F}_p$-point can exist.

Let $Y$ be now a smooth quasiprojective, not necessarily irreducible variety over $\mathbb{F}_p$ and let $Y = Y_1 \cup \cdots \cup Y_r$ be the decomposition into the disjoint union of the irreducible components over $\mathbb{F}_p$. Let them be ordered such that $Y_1, \ldots, Y_l$ are absolutely irreducible and $Y_{l+1}, \ldots, Y_r$ are not. Then the above argument shows that

$$\operatorname{card}(Y(\mathbb{F}_p)) = \operatorname{card}(Y_1(\mathbb{F}_p)) + \cdots + \operatorname{card}(Y_l(\mathbb{F}_p)).$$

The statement of the lemma follows for all primes $p$ such that reduction mod $p$ of $X$ is smooth by the Lang-Weil estimate (see Proposition 4.4). Note the irreducible components of $\overline{X}$ all have dimension $d$, the dimension of $X$, because the irreducible components over $\widehat{\mathbb{F}_p}$ can be obtained from the irreducible components of $X$ over $\widehat{\mathbb{Q}}$ by reduction mod $\mathfrak{p}$ as explained above.

The fact that $l_p(X) > 0$ for a dense set of primes follows from the proof of the previous lemma.                                                                      $\square$

We can use Lemma 4.7 to estimate the size of $c_{p,I}$. Recall that this is the number of points mod $p$ on $E_I \setminus \bigcup_{j \in T \setminus I} E_j$ where $E_I = \bigcap_{i \in I} E_i$.

Let $F_{I,k}, k \in C_I$, be the irreducible components over $\mathbb{Q}$ of $E_I$ with maximal dimension $d_I$ say. We show later (see Proposition 4.13) that $d_I = m - |I|$.

LEMMA 4.8.  *The dimension of every irreducible component of $F_{I,k} \cap \bigcup_{j \in T \setminus I} E_j$ is strictly smaller than $d_I$. (If $F_j$ has dimension zero then this means that $F_I \cap \bigcup_{j \in T \setminus I} E_j$ is the empty set, which by convention has dimension $-\infty$.)*

*Proof.* This follows since the $E_i$ have normal crossings.             $\square$

Since the dimension of $E_I \cap \bigcup_{j \in T \setminus I} E_j$ is strictly smaller than $d_I$, this implies that there exists some $\delta' > 0$ such that for almost all primes $p$

$$\left| c_{p,I} - \sum_{k \in C_I} \operatorname{card} \overline{F_{I,k}}(\mathbb{F}_p) \right| \le \delta' p^{d_I - 1}.$$

But now this together with the Lang-Weil estimate in Lemma 4.7 implies:

PROPOSITION 4.9.  *There exists $\delta \in \mathbb{R}$ such that for almost all primes $p$*

$$|c_{p,I} - \sum_{k \in C_I} l_p(F_{I,k}) p^{d_I}| \le \delta p^{d_I - 1/2}$$

*and $l_p(F_{I,k}) > 0$ for a dense set of primes $p$.*

We make the following definitions for $k \in W'$ and $p \notin Q$:

$$\begin{aligned}
Z_{k,p}(s) &= \frac{c_{p,k}}{a_{p,0}} p^{-m} (p-1)^{|I_k|} \prod_{j \in M_k} \frac{p^{-(A_j s + B_j)}}{1 - p^{-(A_j s + B_j)}}, \\
d_k &= d_I \text{ and } I_k = I \text{ if } k \in W_I.
\end{aligned}$$

If $d_k \ge 0$ then put

$$\alpha_k = \max\left\{\frac{1 + d_k - m + |I_k| - \sum_{j \in M_k} B_j}{\sum_{j \in M_k} A_j}, \frac{-B_j}{A_j} \text{ for } j \in M_k \text{ and } A_j \neq 0\right\}.$$

If $d_k = -\infty$, i.e. $c_{p,k} = 0$ for almost all $p$, then put $\alpha_k = -\infty$.

To analyse the convergence of our Euler product we use the following basic facts:

(A) An infinite product $\prod_{n \in J}(1 + a_n)$ converges absolutely if and only if the corresponding sum $\sum_{n \in J} |a_n|$ converges.

(B) $\sum_p$ prime $|p^{-s}|$ converges if and only if $\Re(s) > 1$. In fact if $\mathcal{P}$ is a dense set of primes we also have that $\sum_{p \in \mathcal{P}} |p^{-s}|$ converges if and only if $\Re(s) > 1$.

(C) The Lang-Weil estimate in Proposition 4.9 means that a sum of the form

$$\sum_p c_{p,I} r_p$$

converges absolutely if and only

$$\sum_p p^{d_I} r_p$$

converges absolutely. Note that we are using here the fact that $l_p(F_{I,k})$ is nonzero with positive density.

The following lemma collects some relevant information about the constant $a_{p,0}$ with which we are normalising our integrals:

LEMMA 4.10. (1) $a_{p,0} = \lim_{s \to \infty} Z_{\mathcal{D}}(s, p)$.

(2)
$$a_{p,0} = \int_{\{\mathbf{x} \in \mathbb{Z}_p^m : \psi(\mathbf{x}) \text{ is valid and } f_0(\mathbf{x}) \in \mathbb{Z}_p^*\}} |g_0(\mathbf{x})| \, |dx|.$$

(3) For each $I \subset T$, $W_I'$ denotes the set of those $k$ for which $\sum_{j \in M_k} A_j \neq 0$. Then for almost all primes $p$

$$a_{p,0} = \sum_{I \subset T} p^{-m} c_{p,I} \sum_{k \in W_I \setminus W_I'} (p-1)^{|I|} \prod_{j \in M_k} \frac{p^{-B_j}}{1 - p^{-B_j}}.$$

Let $T_0$ be the subset of $T$ consisting of those $i$ for which $N_i(f_0) = 0$. Then $W_I' = \emptyset$ for $I \subset T_0$.

(4) There is an integer $N \geq 1$ and a constant $C \geq 0$ such that for all primes $p$

$$1 \geq a_{p,0} \geq 1 - Np^{-1} - Cp^{-3/2}.$$

(5) There exists $M \in \mathbb{N}$ such that for almost all primes $p$

$$1 \leq a_{p,0}^{-1} \leq 1 + Mp^{-1/2}.$$

*Proof.* Most of this is obvious. For (4) note that for almost all $p$, $|g_0(\mathbf{x})|$ $\leq 1$ on $\mathbb{Z}_p^m$ and

$$\mathbb{Z}_p^m \supseteq \left\{ \mathbf{x} \in \mathbb{Z}_p^m : \psi(\mathbf{x}) \text{ holds in } \mathbb{Z}_p \text{ and } f_0(\mathbf{x}) \in \mathbb{Z}_p^* \right\}.$$

Hence from the description of $a_{p,0}$ in (2), $a_{p,0} \leq 1$.

Let $r(\mathbf{x}) := g_0(\mathbf{x})f_0(\mathbf{x})f_1(\mathbf{x}) \cdots f_l(\mathbf{x})$ and

$$L := \{\mathbf{x} \in \mathbb{Z}_p^m : r(\mathbf{x}) \in p\mathbb{Z}_p\}.$$

Note that $\psi(\mathbf{x})$ and $f_0(\mathbf{x}) \in \mathbb{Z}_p^*$ are valid on the complement of $L$ in $\mathbb{Z}_p^m$ and also that $|g_0(\mathbf{x})| = 1$ on this complement.

The equation $r(\mathbf{x}) = 0$ defines a $\mathbb{Q}$-defined hypersurface in $m$-dimensional affine space unless it is the empty set in which case the estimate is trivially valid. There now is an integer $N \geq 1$ so that the number of irreducible components of the reduction of this hypersurface modulo $p$ which are absolutely irreducible is $\leq N$. From the Lang-Weil theorem we know that the complement of $L$ in $\mathbb{Z}_p^m$ contains $\geq p^m - Np^{m-1} - Cp^{m-3/2}$ cosets modulo $p\mathbb{Z}_p$. Since each coset has volume $p^{-m}$ and since $|g_0(\mathbf{x})| = 1$ on this complement the estimate follows from the formula for $a_{p,0}$ in (2). (Of course there is a difficulty in defining the reduction modulo $p$ if the denominators of the $f_i$ and $g_i$ are not prime to $p$. We have to exclude finitely many primes $p$, but the $N$ and $C$ can be chosen to ensure that the estimate holds for all $p$.)

Finally for (5), the estimate in (4) implies that there is a constant $T$ such that $a_{p,0} \geq 1 - Tp^{-1}$ for all but finitely many primes $p$. Hence for almost all $p$:

$$a_{p,0}^{-1} \leq \frac{1}{1 - Tp^{-1}} = 1 + p^{-1/2}\frac{Tp^{-1/2}}{1 - Tp^{-1}}.$$

Since $\frac{Tp^{-1/2}}{1 - Tp^{-1}}$ goes to 0 as $p$ goes to infinity the estimate follows. $\qquad\square$

LEMMA 4.11. *For $k \in W'$, $\alpha_k$ is the abscissa of convergence of*

$$\prod_{p \notin Q} \left(1 + Z_{k,p}(s)\right).$$

*Proof.* First note that each term $\frac{p^{-(A_j s + B_j)}}{1 - p^{-(A_j s + B_j)}}$ converges absolutely if and only if $\Re(s) > \frac{-B_j}{A_j}$.

Suppose that for some $l \in M_k$ for which $A_l \neq 0$,

$$\frac{-B_l}{A_l} \geq \frac{1 + d_k - m + |I_k| - \sum_{j \in M_k} B_j}{\sum_{j \in M_k} A_j}.$$

To prove that $\alpha_k$ is the abscissa of convergence it will suffice in this case to show that

$$(4.3) \qquad \sum_{p \notin Q} (p-1)^{|I_k|} \frac{c_{p,k}}{a_{p,0}} p^{-m} \prod_{j \in M_k} \frac{p^{-(A_j s + B_j)}}{1 - p^{-(A_j s + B_j)}}$$

converges absolutely for $\Re(s) > \alpha_k$. Since $\left(1 - p^{-(A_j s + B_j)}\right)^{-1}$ is a positive decreasing sequence as $p$ increases for $\Re(s) > \frac{-B_j}{A_j}$ if $A_j \neq 0$, it suffices to show that

$$(4.4) \qquad \sum_{p \notin Q} (p-1)^{|I_k|} \frac{c_{p,k}}{a_{p,0}} p^{-m} p^{-\left(\sum_{j \in M_k} A_j s + B_j\right)}$$

converges absolutely.

From Lemma 4.10 (5), $a_{p,0}^{-1} \leq 1 + M p^{-1/2}$. Hence (4.4) converges absolutely provided

$$\sum_{p \notin Q} (p-1)^{|I_k|} c_{p,k} p^{-m} p^{-\left(\sum_{j \in M_k} A_j s + B_j\right)}$$

converges. This sum converges absolutely (by the Lang-Weil argument (C)) if

$$\Re(s) > \frac{1 + d_k - m + |I_k| - \sum_{j \in M_k} B_j}{\sum_{j \in M_k} A_j}.$$

This is the case since $\Re(s) > \alpha_k \geq \frac{1 + d_k - m + |I_k| - \sum_{j \in M_k} B_j}{\sum_{j \in M_k} A_j}$.

Suppose now that $\alpha_k = \frac{1 + d_k - m + |I_k| - \sum_{j \in M_k} B_j}{\sum_{j \in M_k} A_j} > \frac{-B_j}{A_j}$ for all $j \in M_k$.

Now there exist $\varepsilon > 0$, $\delta_1, \delta_2 > 0$ such that $\alpha_k - \varepsilon > \frac{-B_j}{A_j}$ for all $j \in M_k$, and for $\Re(s) > \alpha_k - \varepsilon$ and all primes $p$:

$$1 < \prod_{j \in M_k} \left(1 - p^{-(A_{k,I,j} s + B_{k,I,j})}\right)^{-1} < 1 + \delta_1 p^{-\delta_2}.$$

Hence, when we combine this with Lemma 4.10 (5) and our Lang-Weil argument (C), for $\Re(s) > \alpha_k - \varepsilon$ the sum (4.3) converges absolutely if and only if the sum

$$\sum_{p \notin Q} (p-1)^{|I_k|} p^{d_k} p^{-m} p^{-\left(\sum_{j \in M_k} A_j s + B_j\right)}$$

converges absolutely, i.e., if and only if $\Re(s) > \alpha_k$. This completes the proof that $\alpha_k$ is the abscissa of convergence of $\prod_{p \notin Q} (1 + Z_{k,p}(s))$. $\qquad \square$

COROLLARY 4.12. *The abscissa of convergence of $Z_{\mathcal{D}}(s)$, $\alpha_{\mathcal{D}}$, is equal to*

$$\max \left\{ \{\alpha_k : k \in W'\} \cup \{\beta_p : p \in Q_1\} \right\}$$

*where $\beta_p \in \mathbb{Q}$ is the abscissa of convergence of the exceptional local factor $P_p(p^{-s})$ where $Q_1$ is the set of bad primes. In particular, $\alpha_{\mathcal{D}}$ is a rational number.*

We shall improve this description of $\alpha_{\mathcal{D}}$ in Lemma 4.15.

The next task is to show that we can meromorphically continue $Z_{\mathcal{D}}(s)$ a little beyond $\Re(s) > \alpha_{\mathcal{D}}$ to allow an application of the Tauberian theorem. The first step is to show that the maximal value of $\alpha_k$ as $k$ varies over the set $W' \subset W$ (where recall $W$ indexes the open simplicial pieces $R_k$ of $\overline{D_T}$ ) is realised only by the one-dimensional simplicial edges, i.e. $k \in \{1, \ldots, q\}$ .

The essential fact that we shall use here is the following:

PROPOSITION 4.13.    *For $k \in W$, $d_k = m - |I_k|$.*

*Proof.* Each $E_i$ is of dimension $m - 1$. Since the $E_i$ intersect with normal crossings, the irreducible components of $\bigcap_{i \in I_k} E_i$ have dimension $d_k = m - |I_k|$. For further details see the proof of Theorem 2.4 of [5] or 17.F of [33].    □

This already gives a much simpler description of $\alpha_k$ :

COROLLARY 4.14.    (1) *For $k \in W'$,*

$$\alpha_k = \max\left\{ \frac{1 - \sum_{j \in M_k} B_j}{\sum_{j \in M_k} A_j}, \frac{-B_j}{A_j} \text{ for } j \in M_k \text{ and } A_j \neq 0 \right\}.$$

(2) *If $k \in W' \cap \{1, \ldots, q\}$ then*

$$\alpha_k = \frac{1 - B_k}{A_k}.$$

LEMMA 4.15.

$$\alpha_{\mathcal{D}} = \max\left\{ \alpha_k : k \in W' \right\} = \max\left\{ \alpha_k : k \in W' \cap \{1, \ldots, q\} \right\}$$

*and $\alpha_k < \alpha_{\mathcal{D}}$ if $k > q$.*

*Proof.* The first fact to point out is that since the abscissa of convergence of the local factors $\beta_p$ for the bad primes $p$ is one of $-B_k/A_k$ for $k \in W' \cap \{1, \ldots, q\}$,

$$\begin{aligned} \alpha_{\mathcal{D}} &= \max\left\{ \{\alpha_k : k \in W'\} \cup \{\beta_p : p \in Q_1\} \right\} \\ &= \max\left\{ \alpha_k : k \in W' \right\} \end{aligned}$$

by (2) of the previous corollary. The rest follows once we can show that for $k > q$ and $k \in W'$

$$(4.5) \qquad \frac{1 - \sum_{j \in M_k} B_j}{\sum_{j \in M_k} A_j} < \max\left\{ \frac{1 - B_j}{A_j} : j \in M_k \cap W' \right\}.$$

Recall the definition of $W'$ from Lemma 4.10 that $k \in W'$ if $\sum_{j \in M_k} A_j \neq 0$. We note that if $j \in M_k$ then $j \in W'$ since

$$\frac{1 - \sum_{j \in M_k} B_j}{\sum_{j \in M_k} A_j} \leq \frac{1 - \sum_{j \in M_k \cap W'} B_j}{\sum_{j \in M_k \cap W'} A_j}.$$

We may suppose without loss of generality that $|M_k| = 2$, $\frac{1-B_2}{A_2} \leq \frac{1-B_1}{A_1}$ and that $A_1$ and $A_2 > 0$. Then

$$\frac{1-(B_1+B_2)}{A_1+A_2} < \frac{1-B_1}{A_1}$$

if and only if

$$A_1 - A_1 B_1 - A_1 B_2 < A_1 + A_2 - A_1 B_1 - A_2 B_1.$$

But our assumptions that $\frac{1-B_2}{A_2} \leq \frac{1-B_1}{A_1}$ and $A_2 > 0$ imply this second inequality. This confirms statement (4.5) and hence the lemma. $\square$

THEOREM 4.16. $\prod_{p \notin Q} \left( 1 + \sum_{j=1}^{w} Z_{j,p}(s) \right)$ has a meromorphic continuation to $\Re(s) > \alpha_\mathcal{D} - \delta$ for some $\delta > 0$.

*Proof.* Define

$$R = \left\{ k \in W' \cap \{1, \ldots, q\} : \frac{1-B_j}{A_j} = \alpha_\mathcal{D} \right\},$$

$$V_p(s) = \prod_{k \in R} \left( 1 - c_{p,k} p^{-m+|I_k|} p^{-(A_k s + B_k)} \right).$$

For convenience, define $Z_{j,p}(s) = 0$ for the finite number of primes $p \in Q$. Hence $\prod_{p \notin Q} \left( 1 + \sum_{j=1}^{w} Z_{j,p}(s) \right) = \prod_p \left( 1 + \sum_{j=1}^{w} Z_{j,p}(s) \right)$.

We introduce the following notation which will be convenient during the course of the proof. Write $\prod_p F_p(s) \equiv \prod_p G_p(s)$ if there exists $\delta > 0$ such that $\sum_p (F_p(s) - G_p(s))$ converges for $\Re(s) > \alpha_\mathcal{D} - \delta$. To prove the lemma it will suffice to prove the following:

(1) $\prod_p V_p(s)$ is a meromorphic function on $\Re(s) > \alpha_\mathcal{D} - \delta$ for some $\delta > 0$.

(2) $\prod_p \left( 1 + \sum_{j=1}^{w} Z_{j,p}(s) \right) V_p(s) \equiv 1$; i.e., the Euler product converges on $\Re(s) > \alpha_\mathcal{D} - \delta$ for some $\delta > 0$.

We prove (2). It will be convenient to note the following fact: suppose that for $s > \alpha_\mathcal{D} - \delta$, the function $X_p(s)$ as $p \to \infty$ is a positive decreasing sequence. If $\prod_p F_p(s) \equiv \prod_p G_p(s)$ then $\prod_p F_p(s) X_p(s) \equiv \prod_p G_p(s) X_p(s)$.

By Lemma 4.15, $\prod_p \left( 1 + \sum_{j=1}^{w} Z_{j,p}(s) \right) \equiv \prod_p \left( 1 + \sum_{k \in R} Z_{k,p}(s) \right)$. Recall the definition of $Z_{k,p}(s)$ for $k \in \{1, \ldots, q\}$ and $p \notin Q$:

$$Z_{k,p}(s) = \frac{c_{p,k}}{a_{p,0}} p^{-m} (p-1)^{|I_k|} \frac{p^{-(A_k s + B_k)}}{1 - p^{-(A_k s + B_k)}}.$$

Note that $\sum_p c_{p,k} p^{-m+|I_k|} \frac{p^{-2(A_k s + B_k)}}{1 - p^{-(A_k s + B_k)}}$ converges on $\Re(s) > \alpha_\mathcal{D} - \delta$ for some $\delta > 0$. This follows because $\left( 1 - p^{-(A_k s + B_k)} \right)^{-1}$ is a positive decreasing sequence and $\frac{(1-2B_k)}{2A_k} < \frac{1-B_k}{A_k}$. Putting this together with our analysis on the estimate for $a_{p,0}$ we see that $\prod_p \left(1 + \sum_{k \in R} Z_k(s)\right) \equiv \prod_p \left(1 + \sum_{k \in R} \overline{Z_k(s)}\right)$ where, for all $p$ now,

$$\overline{Z_{k,p}(s)} = c_{p,k} p^{-m+|I_k|} p^{-(A_k s + B_k)}.$$

Hence we have

$$(4.6) \qquad \prod_p \left(1 + \sum_{j=1}^{w} Z_{j,p}(s)\right) \equiv \prod_p \left(1 + \sum_{k \in R} \overline{Z_{k,p}(s)}\right).$$

Now

$$(4.7)$$
$$\prod_p V_p(s) \equiv \prod_p \left(1 - \sum_{k \in R} c_{p,k} p^{-m+|I_k|} p^{-(A_k s + B_k)}\right) = \prod_p \left(1 - \sum_{k \in R} \overline{Z_{k,p}(s)}\right)$$

since $\frac{1 - B_{k_1} - \cdots - B_{k_i}}{A_{k_1} + \cdots + A_{k_i}} < \frac{1 - B_k}{A_k}$ for $i > 1$ and $k_j, k \in R$. Now from (4.6) and (4.7) we obtain

$$\prod_p \left(1 + \sum_{j=1}^{w} Z_{j,p}(s)\right) V_p(s) \equiv \prod_p \left(1 + \sum_{j=1}^{w} Z_{j,p}(s)\right) \left(1 - \sum_{k \in R} \overline{Z_{k,p}(s)}\right)$$
$$\equiv \prod_p \left(1 + \sum_{k \in R} \overline{Z_{k,p}(s)}\right) \left(1 - \sum_{k \in R} \overline{Z_{k,p}(s)}\right).$$

Here we use the fact that $\sum_{j=1}^{w} Z_{j,p}(s)$ and $\sum_{k \in R} c_{p,k} p^{-m+|I_k|} p^{-(A_k s + B_k)}$ are both decreasing positive sequences as $p \to \infty$. Finally the fact that $\frac{(1 - B_{k_1} - B_{k_2})}{A_{k_1} + A_{k_2}} < \frac{1 - B_{k_i}}{A_{k_i}}$ implies that

$$\prod_p \left(1 + \sum_{k \in R} \overline{Z_{k,p}(s)}\right) \left(1 - \sum_{k \in R} \overline{Z_{k,p}(s)}\right) \equiv 1.$$

Hence we have proved (2).

For the proof of (1) note that by Proposition 4.9

$$\prod_{p \text{ prime}} \left(1 - c_{p,k} p^{-m+|I_k|} p^{-(A_k s + B_k)}\right)$$

$$\equiv \prod_{p \text{ prime}} \left(1 - \sum_{j \in C_I} l_p(F_{I_k,j}) p^{d_k} p^{-m+|I_k|} p^{-(A_k s + B_k)}\right)$$

$$\equiv \prod_{j \in C_I} \prod_{p \text{ prime}} \left(1 - l_p(F_{I_k,j}) p^{-(A_k s + B_k)}\right).$$

Hence the fact that $\prod_p V_p(s)$ is meromorphic on $\Re(s) > \alpha_\mathcal{D} - \delta$ for some $\delta > 0$ follows from Lemma 4.6. This completes the proof of Theorem 4.16. $\qquad\square$

To prove our main theorem about the growth of the coefficients in a Dirichlet series defined as the Euler product of cone integrals over $\mathbb{Q}$, we can now apply the Tauberian theorems. First of all, an easy consequence of the Hardy-Littlewood-Karamata Tauberian theorem as formulated in [38, §7.3, Th. 5] is the following:

THEOREM 4.17.  *Let $F(s) = \sum_{n=1}^{\infty} a_n n^{-s}$ be a Dirichlet series convergent for $s \in \mathbb{C}$ with $\Re(s) > \alpha \geq 0$. Suppose that $a_n \in \mathbb{R}$ and $a_n \geq 0$ for all $n \in \mathbb{N}$. Suppose further that*

$$F(\sigma) = (c + o(1)) (\sigma - \alpha)^{-w}$$

*for $\sigma \in \mathbb{R}$, $\sigma > \alpha$ and $\sigma \to \alpha$. Then,*

$$\sum_{n=1}^{N} \frac{a_n}{n^{\alpha}} \sim \frac{c}{\Gamma(w+1)} \cdot (\log N)^w$$

*for $N \to \infty$.*

This theorem implies:

COROLLARY 4.18.  *Let $Z(s) = \sum_{n=1}^{\infty} a_n n^{-s}$ be defined as an Euler product of cone integrals over $\mathbb{Q}$. Suppose $Z(s)$ is not the constant function.*

*(1) The abscissa of convergence $\alpha$ of $Z(s)$ is a rational number and $Z(s)$ has a meromorphic continuation to $\Re(s) > \alpha - \delta$ for some $\delta > 0$.*

*(2) Let the pole at $s = \alpha$ have order $w$. Then,*

$$\sum_{n=1}^{N} \frac{a_n}{n^{\alpha}} \sim \frac{c}{\Gamma(w+1)} \cdot (\log N)^w$$

*for $N \to \infty$ and some real number $c \in \mathbb{R}$.*

Another Tauberian theorem that we can apply is Ikehara's theorem as formulated on page 62 of [2]:

THEOREM 4.19.  *Suppose $F(s) = \frac{1}{s} \int_0^{\infty} e^{-st} dA(t)$ for $\Re(s) > a > 0$, $A(t) > 0$. Suppose further that in a neighbourhood of $s = a$*

$$F(s) = g(s)(s-a)^{-w} + h(s)$$

*where $g$ and $h$ are holomorphic and $g(a) \neq 0$. Assume also that $F(s)$ can be holomorphically continued to the line $\Re(s) = a$ except for the pole at $s = a$. Then*

$$A(t) \sim \frac{g(a)}{\Gamma(w)} \cdot e^{at} t^{w-1}.$$

In our case, $A(t) = \sum_{\log n < t} a_n$ and $f(s) = \sum_{n=1}^{\infty} a_n n^{-s} = sF(s)$; then substituting $t = \log x$ we get:

THEOREM 4.20.  *Let the Dirichlet series $f(s) = \sum_{n=1}^{\infty} a_n n^{-s}$ with non-negative coefficients be convergent for $\Re(s) > a > 0$. Assume in a neighbourhood of $a$, $f(s) = g(s)(s-a)^{-w} + h(s)$, holds where $g(s), h(s)$ are holomorphic functions, $g(a) \neq 0$ and $w > 0$. Assume also that $f(s)$ can be holomorphically continued to the line $\Re(s) = a$ except for the pole at $s = a$. Then for $x$ tending to infinity, we have*

$$\sum_{n \leq x} a_n \sim \left( \frac{g(a)}{a\Gamma(w)} \right) \cdot x^a (\log x)^{w-1}.$$

To apply this Tauberian theorem we shall need the following facts about the Artin $L$-function:

PROPOSITION 4.21.  *Let $\mathbb{G}$ be the absolute Galois group of $\mathbb{Q}$ and $\rho : \mathbb{G} \to \mathrm{GL}(V)$ be a continuous finite-dimensional complex representation of $\mathbb{G}$. Let*

$$L(\rho, s) := \prod_p \frac{1}{\det_{V^{I_p}}(1 - \rho(\mathrm{Frob}_p) \cdot p^{-s})}$$

*be the corresponding Artin $L$-function defined as an Euler product over all primes $p$ (see [32]). Then the following hold:*

(1)  *The poles of the Euler factors of $L(\rho, s)$ are on the line $\Re(s) = 0$.*

(2)  *The Dirichlet series $L(\rho, s)$ converges for $\Re(s) > 1$ and has meromorphic continuation to all of $\mathbb{C}$.*

(3)  *The extension of $L(\rho, s)$ has no pole or zero on the line $\Re(s) = 1$ except possibly in $s = 1$.*

*Proof.*  Properties (1) and (2) are well known; see [32]. To prove (3) note that by the Brauer induction theorem (see [32]) and by class field theory, there are Hecke characters $\chi_1, \ldots, \chi_k$ and $\psi_1, \ldots, \psi_l$ of appropriate number fields such that

$$L(\rho, s) = \frac{L(\chi_1, s) \cdots L(\chi_k, s)}{L(\psi_1, s) \cdots L(\psi_l, s)}.$$

The Hecke $L$-functions appearing in this formula have the desired property (3) (see [28, Ch. 15, §4] and [3, Ch. 13] for the case of trivial characters).  □

We deduce from this proposition something about our Euler products.

COROLLARY 4.22.  *Let $Z_{\mathcal{D}}(s) = \sum_{n=1}^{\infty} a_n n^{-s}$ be defined as an Euler product*

$$Z_{\mathcal{D}}(s) = \prod_{p \ prime, \ a_{p,0} \neq 0} \left( a_{p,0}^{-1} \cdot Z_{\mathcal{D}}(s, p) \right)$$

*of cone integrals over $\mathbb{Q}$. Suppose $Z_{\mathcal{D}}(s)$ is not the constant function.*

(1) *The abscissa of convergence $\alpha = \alpha_{\mathcal{D}}$ of $Z_{\mathcal{D}}(s)$ is a rational number and $Z_{\mathcal{D}}(s)$ can be holomorphically continued to the line $\Re(s) = \alpha$ except for the pole at $s = \alpha$.*

(2) *Let the pole at $s = \alpha$ have order $w$. Then there exists some real number $c \in \mathbb{R}$ such that*

$$(4.8) \qquad\qquad a_1 + a_2 + \cdots + a_N \sim c \cdot N^{\alpha} \left( \log N \right)^{w-1}$$

*as $N \to \infty$.*

*Proof.* The proof of Theorem 4.16 established that $Z_{\mathcal{D}}(s)$ is the product of Artin $L$-functions and a Dirichlet series convergent for $\Re(s) > \alpha_{\mathcal{D}} - \delta$ for some $\delta > 0$. Note that for those primes $p \in Q_1$ with bad reduction, our explicit expression (3.2) implies that $\beta_p < \alpha_{\mathcal{D}}$ where $\beta_p$ was the abscissa of convergence of $Z_{\mathcal{D}}(s, p)$. $\qquad\square$

It was important to establish that the abscissa of convergence of the local factors is strictly to the left of $\alpha_{\mathcal{D}}$. For example the coefficients of $\zeta(s) \cdot (1 - p^{s-1})$ do not satisfy an asymptotic formula of the form (4.8). This was pointed out to us by Benjamin Klopsch and Dan Segal.

Note that we do not know anything about the possibility of other poles in the region $\alpha = \Re(s) > \alpha - \delta$ other than the one at $s = \alpha$. This is because we have used Artin $L$-functions to do our meromorphic continuation. However it is conjectured that these Artin $L$-functions actually have this one pole and no others (see [32]).

## 5. Nilpotent groups

In this section we show that the zeta functions of finitely generated nilpotent groups are essentially Euler products of cone integrals over $\mathbb{Q}$.

To simplify our analysis of subgroups and normal subgroups we use the following notation: $\zeta_G^{\leq}(s) = \zeta_G(s)$.

The following extends a result in [24]. It will allow us to concentrate on counting in rings rather than in groups:

**THEOREM 5.1.** *Let $G$ be a finitely generated nilpotent group. Then there is a subgroup of finite index $G_0$ and a Lie algebra $L(G_0)$ over $\mathbb{Z}$ constructed as the image under $\log$ of $G_0$ such that:*

(1) *For $* \in \{\leq, \lhd\}$ and almost all primes $p$*

$$\zeta_{G,p}^*(s) = \zeta_{G_0,p}^*(s) = \zeta_{L(G_0),p}^*(s).$$

(2) *If $\alpha_p^*(G), \alpha_p^*(G_0)$ and $\alpha_p^*(L(G_0))$ denote the abscissas of convergence of the local factors $\zeta_{G,p}^*(s), \zeta_{G_0,p}^*(s)$ and $\zeta_{L(G_0),p}^*(s)$ respectively, then for all*

*primes p*

$$\alpha_p^*(G) = \alpha_p^*(G_0) = \alpha_p^*(L(G_0)).$$

*Proof.* Part (1) follows from Section 4 of [24]. The same paper contains a proof (of Proposition 1.8) that the abscissa of convergence is a commensurability invariant which implies that $\alpha_p^*(G) = \alpha_p^*(G_0)$. (Note that the statement of Proposition 1.8 is slightly weaker applying to torsion-free nilpotent groups, but the proof suffices to establish that $\alpha_p^*(G) = \alpha_p^*(G_0)$.) A similar (slightly simpler) argument shows that if $L_0$ is a Lie subring of finite index in $L$ then $\alpha_p^*(L) = \alpha_p^*(L_0)$.

We shall deduce the remaining equality $(\alpha_p^*(G_0) = \alpha_p^*(L(G_0)))$ from Proposition 5.2 below. We shall explain this argument together with the construction of $G_0$ from $G$. We consider the prime $p$ to be fixed and write $\widehat{G}_p$ for the pro-$p$ completion of $G$. We have $\alpha_p^*(G) = \alpha_p^*(\widehat{G}_p)$ and $\alpha_p^*(L) = \alpha_p^*(L \otimes \mathbb{Z}_p)$ ($L$ a Lie algebra). The abscissas of convergence $\alpha_p^*(\widehat{G}_p)$, $\alpha_p^*(L \otimes \mathbb{Z}_p)$ are those of the power series (in $p^{-s}$) counting open subgroups, open subalgebras of finite index respectively or open normal subgroups and open ideals of finite index. We may now replace $G$ by $G_0$ chosen so that $\widehat{G}_{0p}$ is a uniform pro-$p$ group (see [7]). The result now follows by a straightforward application of the proposition below. □

The following proposition is due to Dan Segal. We cordially thank him for the permission to include its proof here. Let $G$ be a pro-$p$ group and $L$ be a finitely generated Lie algebra over $\mathbb{Z}_p$. We write in accordance with previous notations $s_n^{\leq}(G)$, $s_n^{\leq}(L)$ for the number of open subgroups, open subalgebras of index at most $p^n$ in $G$, $L$ respectively, and $s_n^{\triangleleft}(G)$, $s_n^{\triangleleft}(L)$ for the corresponding numbers of open normal subgroups of $G$, open ideals of $L$. We shall use without special mention the results of [7, §9.4], concerning the correspondence between powerful Lie algebras and uniform pro-$p$ groups. Note in particular that $L = \log G$ is a powerful $\mathbb{Z}_p$-Lie algebra for any nilpotent uniform pro-$p$ group $G$ (see [7, §9.4]).

PROPOSITION 5.2. *Let $G$ be a nilpotent uniform pro-$p$ group and define $L = \log G$ to be its Lie algebra. There exist $a, b, B, D \in \mathbb{N}$, an open normal subgroup $G_1$ of $G$ and an open powerful subalgebra $L_1$ of $L$ such that, for all $n$,*

$$
\begin{align}
(5.1) && s_n^{\leq}(L) &\leq s_{n+a}^{\leq}(G) \\
(5.2) && s_n^{\leq}(G) &\leq B s_{n+b}^{\leq}(L) \\
(5.3) && s_n^{\triangleleft}(L) &\leq s_n^{\triangleleft}(G_1) \\
(5.4) && s_n^{\triangleleft}(G) &\leq D s_n^{\triangleleft}(L_1).
\end{align}
$$

*Proof.* Say $G$ has nilpotency class $c$ and $\dim L = r$. We put $\mathbf{p} = p$ if $p$ is odd, $\mathbf{p} = 4$ if $p = 2$.

(5.1). If $A$ is an open subalgebra of $L$ then $\mathbf{p}A$ is powerful so $H(A) := \exp(\mathbf{p}A)$ is a uniform subgroup of $G$, and

$$|G : H(A)| = |L : \mathbf{p}A| = \mathbf{p}^r |L : A|.$$

Since $A \mapsto H(A)$ is one-to-one it follows that $s_{\bar{n}}^{\leq}(L) \leq s_{n+a}^{\leq}(G)$ where $a = r$ if $p$ is odd, $a = 2r$ if $p = 2$. (This part does not depend on $G$ being nilpotent.)

(5.2). There exists $f \in \mathbb{N}$ such that for each open subgroup $H$ of $G$, the subgroup $H^{p^f}$ is uniform ([7, Prop. 3.9 and Thm. 4.5]). Then $M(H) := \log(H^{p^f})$ is a powerful subalgebra of $L$, and

$$|L : M(H)| = \left|G : H^{p^f}\right| \leq p^b |G : H|$$

where $b = fr$. The mapping $H \mapsto M(H)$ need not be one-to-one; we show that its fibres have bounded order. The subgroup $H$ is contained in

$$H_1 = \left\langle \exp(p^{-f} M(H)) \right\rangle,$$

and $H_1$ is generated by elements $x$ such that $x^{p^f} \in H^{p^f}$, which implies that $H_1^{p^m} \leq H^{p^f} \leq H$ where $m = c(c+1)f/2$ ([35, Chap. 6, Prop. 3]). The group $H_1/H_1^{p^m}$ has order at most $p^{mr}$ and rank at most $r$, hence contains at most $p^{mr^2}$ subgroups. Hence the number of subgroups $H$ corresponding to a given $M(H)$ is at most $p^{mr^2}$, and it follows that $s_{\bar{n}}^{\leq}(G) \leq B s_{n+b}^{\leq}(L)$ where $B = p^{mr^2}$.

(5.3). This depends on the 'commutator Campbell-Hausdorff formula' (see [7, §6.3]): for $h, x \in G$ with $\log h = a$, $\log x = u$ we have

$$(5.5) \qquad \log[h, x] = (a, u) + \sum q_{\mathbf{e}}(a, u)_{\mathbf{e}},$$

where each $(a, u)_{\mathbf{e}}$ is a repeated Lie bracket involving at least one $a$ and one $u$, the $q_{\mathbf{e}}$ are rational numbers depending only on $\mathbf{e}$, and the sum is finite since $L$ is nilpotent. Fix $f \in \mathbb{N}$ so that $p^f q_{\mathbf{e}} \in \mathbb{Z}$ for each $\mathbf{e}$ (and assume that $f \geq 2$ if $p = 2$). Now put $G_1 = G^{p^f}$.

Let $I$ be an open ideal of $L$. Then $p^f I$ is a powerful Lie subalgebra of $L$, so $H(I) := \exp(p^f I)$ is a uniform subgroup of $G$, contained in $G_1$. Moreover, applying (5.5) with $h \in H(I)$ and $x = y^{p^f} \in G_1$ we see that

$$\log[h, x] = (a, u) + \sum p^{n(\mathbf{e})f} q_{\mathbf{e}}(a, \log y)_{\mathbf{e}} \in p^f I$$

where $n(\mathbf{e})$ denotes the number of occurrences of $u$ in $(a, u)_{\mathbf{e}}$. It follows that $H(I)$ is normal in $G_1$. Also

$$|G_1 : H(I)| = \left|\log G_1 : p^f I\right| = \left|p^f L : p^f I\right| = |L : I|.$$

As $I \mapsto H(I)$ is one-to-one this shows that $s_n^{\triangleleft}(L) \leq s_n^{\triangleleft}(G_1)$.

(5.4). According to [7, Prop. 3.9], $G$ has an open normal subgroup $G_1$ such that whenever $T$ is an open normal subgroup of $G$ contained in $G_1$, $T$ is powerfully embedded in $G_1$ (i.e. $[T, G_1] \leq T^{\mathbf{p}}$). In particular, $T$ is a uniform group. By a result of Shalev (see [7, Exercise 2.2(iii)]) we then have

$$[T^{p^n}, G_1^{p^n}] \leq [T, G_1]^{p^{2n}} \leq T^{p^{2n}\mathbf{p}}$$

for all $n \in \mathbb{N}$. This implies that for $x \in T$ and $y \in G_1$ we have

$$(\log x, \log y) \in \mathbf{p} \log T,$$

by [7], Lemma 7.12. It follows that $\log T$ is an ideal of $L_1 = \log G_1$.

Now to each open normal subgroup $N$ of $G$ we associate

$$I(N) := \log(N \cap G_1).$$

The preceding paragraph shows that $I(N)$ is an ideal of $L_1$, and we have

$$|L_1 : I(N)| = |G_1 : N \cap G_1| \leq |G : N|.$$

It remains to bound the fibres of the mapping $N \mapsto I(N)$. Given $N \cap G_1 = N_1$ and $NG_1 = N_2$, the number of possibilities for $N$ is at most $|\mathrm{Hom}(N_2/G_1, G_1/N_1)|$. Say $p^k$ is the exponent of $G/G_1$. Since every characteristic subgroup of $G_1/N_1$ is powerful, Exercise 2.5(d) of [7] shows that the elements of order dividing $p^k$ in $G_1/N_1$ form a subgroup $E_k(G_1/N_1) = E$, say; also $E$ is powerful and has rank at most $r$, so $|E| \leq p^{kr}$. Since $G/G_1$ also has rank at most $r$ it follows that

$$|\mathrm{Hom}(N_2/G_1, G_1/N_1)| = |\mathrm{Hom}(N_2/G_1, E)| \leq p^{kr^2}.$$

We conclude that the number of normal subgroups $N$ of $G$ giving rise to a given $I(N)$ is at most $D = D_1 p^{kr^2}$ where $D_1$ is the number of normal subgroups in the finite group $G/G_1$, and it follows that $s_n^\lhd(G) \leq D s_n^\lhd(L_1)$.                □

The finitely many exceptional primes in Theorem 5.1 will not worry us thanks to the following result, proved first for torsion-free nilpotent groups in [24] and more generally by the first author in [8]:

THEOREM 5.3.   *Let $G$ be a finitely generated group of finite rank. Then, for $* \in \{\leq, \lhd\}$ and for each prime $p$, $\zeta_{G,p}^*(s)$ is a rational function.*

Recall that a group $G$ has finite rank if there is a bound on the number of generators of finitely generated subgroups; for example this holds if $G$ is a finitely generated nilpotent group or more generally a polycyclic group.

The proof of the rationality of the local factors depended on expressing them as "definable" $p$-adic integrals. We recall the description of these integrals in the case of a ring $L$ which is additively isomorphic to $\mathbb{Z}^d$. Fix a basis for

$L$ which identifies it with $\mathbb{Z}^d$. The multiplication in $L$ is given by a bilinear mapping

$$\beta : \mathbb{Z}^d \times \mathbb{Z}^d \to \mathbb{Z}^d,$$

which extends to a bilinear map $\mathbb{Z}_p^d \times \mathbb{Z}_p^d \to \mathbb{Z}_p^d$ for each prime $p$, giving the structure of the $\mathbb{Z}_p$-algebra $L_p = L \otimes \mathbb{Z}_p$.

Define

$$V_p^{\leq} = \left\{ \begin{array}{c} (m_{kl}) \in \mathrm{Tr}_d\left(\mathbb{Z}_p\right) : \text{for } 1 \leq i,j \leq d, \exists Y_{ij}^1, \ldots, Y_{ij}^d \in \mathbb{Z}_p \\ \text{such that } \beta(\mathbf{m}_i, \mathbf{m}_j) = \sum_{k=1}^d Y_{ij}^k \mathbf{m}_k \end{array} \right\}$$

and

$$V_p^{\triangleleft} = \left\{ \begin{array}{c} (m_{kl}) \in \mathrm{Tr}_d\left(\mathbb{Z}_p\right) : \text{for } 1 \leq i,j \leq d, \exists Y_{ij}^1, \ldots, Y_{ij}^d \in \mathbb{Z}_p \\ \text{such that } \beta(\mathbf{m}_i, \mathbf{e}_j) = \sum_{k=1}^d Y_{ij}^k \mathbf{m}_k \end{array} \right\}$$

where $\mathrm{Tr}_d\left(\mathbb{Z}_p\right)$ denotes the set of upper triangular matrices. Each subset consists of matrices whose rows $\mathbf{m}_i$ generate a lattice in $\mathbb{Z}_p^d$ which is either a subalgebra or an ideal. The $d$ tuple $\mathbf{e}_j$ denotes the standard unit vector with 1 in the $j^{\text{th}}$ entry and zeros elsewhere and corresponds to the $j^{\text{th}}$ basis element.

PROPOSITION 5.4 ([24, Prop. 3.1]). *For $* \in \{\leq, \triangleleft\}$ and each prime $p$*

$$\zeta_{L,p}^*(s) = (1 - p^{-1})^{-d} \int_{V_p^*} |m_{11}|^{s-1} \cdots |m_{dd}|^{s-d} \, |dx|$$

*where $v$ is the valuation on $\mathbb{Z}_p$, $|m| = p^{-v(m)}$ and $|dx|$ is the normalized additive Haar measure on $\mathbb{Z}_p^{d(d+1)/2} \equiv \mathrm{Tr}_d\left(\mathbb{Z}_p\right)$.*

The proof of the rationality of these $p$-adic integrals relies on observing that $V_p^*$ are definable subsets in the language of fields. One can then apply a theorem of Denef's [4] which establishes the rationality of definable $p$-adic integrals. Denef's proof relies on an application of Macintyre's quantifier elimination for the theory of $\mathbb{Q}_p$ which simplifies in a generally mysterious way the description of definable subsets like $V_p^*$. However in our case it is possible to do this elimination by hand since it involves solving linear equations. This removes the necessity for the model theoretic black box in the proof of the rationality.

THEOREM 5.5. *Let $L$ be a ring additively isomorphic to $R^d$ where $R = \mathbb{Z}$ or $\mathbb{Z}_p$. For $* \in \{\leq, \triangleleft\}$ there exist homogeneous polynomials*

$$g_{ijk}^*(\mathbf{X}) \in R[X_{rs} : 1 \leq r \leq s \leq d]$$

*($i,j,k \in \{1, \ldots, d\}$) of degree $k$ for $* = \triangleleft$ and $k+1$ for $* = \leq$ such that*

$$V_p^* = \Big\{ (m_{rs}) \in \mathrm{Tr}_d\left(\mathbb{Z}_p\right) : v(m_{11} \cdots m_{kk})$$
$$\leq v(g_{ijk}^*(m_{rs})) \text{ for } i,j,k \in \{1, \ldots, d\} \Big\}.$$

*Proof.* We can express the defining conditions for $V_p^*$ in matrix form which makes things quite transparent. Let $C_j$ denote the matrix whose rows are $\mathbf{c}_i = \beta(\mathbf{e}_i, \mathbf{e}_j)$.

$M \in V_p^{\lhd}$ if we can solve for each $1 \le i, j \le d$ the equation

$$\mathbf{m}_i C_j = \left( Y_{ij}^1, \dots, Y_{ij}^d \right) M$$

with $\left( Y_{ij}^1, \dots, Y_{ij}^d \right) \in \mathbb{Z}_p^d$. Let $M'$ denote the adjoint matrix and

$$M^{\natural} = M' \mathrm{diag}(m_{22}^{-1} \cdots m_{dd}^{-1}, \dots, m_{dd}^{-1}, 1).$$

Then since the matrix $M$ is upper triangular, the $ik^{\mathrm{th}}$ entry of $M^{\natural}$ is a homogeneous polynomial of degree $k - 1$ in the variables $m_{rs}$ with $1 \le r \le s \le k - 1$. Then we can rewrite the above equation as:

$$\mathbf{m}_i C_j M^{\natural} = \left( m_{11} Y_{ij}^1, \dots, m_{11} \cdots m_{dd} Y_{ij}^d \right).$$

Let $g_{ijk}^{\lhd}(m_{rs})$ denote the $k^{\mathrm{th}}$ entry of the $d$-tuple $\mathbf{m}_i C_j M^{\natural}$ which is a homogeneous polynomial of degree $k$ in $m_{rs}$. (In fact we can see that it is homogeneous of degree 1 in $m_{is}$ ( $s = 1, \dots, d$) and degree $k - 1$ in $m_{rs}$ with $1 \le r \le s \le k - 1$.) Then $V_p^{\lhd}$ has the description detailed in the statement of the theorem.

For $M \in V_p^{\le}$ we are required to solve for each $1 \le i, j \le d$ the equation

$$\mathbf{m}_i \left( \sum_{l=j}^{d} m_{jl} C_l \right) M^{\natural} = \left( m_{11} Y_{ij}^1, \dots, m_{11} \cdots m_{dd} Y_{ij}^d \right)$$

with $\left( Y_{ij}^1, \dots, Y_{ij}^d \right) \in \mathbb{Z}_p^d$. Let $g_{ijk}^{\le}(m_{rs})$ denote the $k^{\mathrm{th}}$ entry of the $d$-tuple $\mathbf{m}_i \left( \sum_{l=j}^{d} m_{jl} C_l \right) M^{\natural}$ which is a homogeneous polynomial of degree $k + 1$ in $m_{rs}$. Again with these polynomials, $V_p^{\le}$ has the description detailed in the statement of the theorem.  □

Now, our theorem has the following corollary stating that the local zeta functions of the ring $L$ can be represented by cone integrals over $\mathbb{Q}$:

COROLLARY 5.6. *Let $\psi^*$ be the cone condition defining $V_p^*$. Set $f_0 = m_{11} \cdots m_{dd}$ and $g_0 = m_{11}^{d-1} \cdots m_{d-1,d-1}$. Let $\mathcal{D}^*$ be the associated cone integral data. Then*

$$\zeta_{L,p}^*(s) = (1 - p^{-1})^{-d} Z_{\mathcal{D}^*}(s - d, p).$$

We know that the constant term of $Z_{\mathcal{D}^*}(s - d, p)$ must be $(1 - p^{-1})^d$ since $\zeta_{L,p}^*(s)$ has constant term 1. So, by Theorem 5.1 (1) the global zeta function $\zeta_G^*(s)$ is the Euler product of the cone integrals up to multiplication by a finite number of rational functions in $p^{-s}$. Note that the abscissas of convergence of these finite numbers of rational functions are strictly less than the abscissa of convergence of this Euler product by Theorem 5.1 (2) and the observation of the previous section that the abscissa of convergence of each $Z_{\mathcal{D}^*}(s,p)$ is strictly less than the abscissa of convergence $\alpha_{\mathcal{D}^*}$ of $Z_{\mathcal{D}^*}(s)$. This is important in application of the second Tauberian theorem of the previous section since we must guarantee that there are no poles creeping in from the exceptional local factors which lie on $\Re(s) = \alpha_{\mathcal{D}^*}$. Hence we get as a corollary of our work on the zeta functions defined as Euler products of cone integrals the following:

THEOREM 5.7. (1) *Let $G$ be a finitely generated nilpotent infinite group. Then there exist finitely many varieties $E_i$, $i \in T$, defined over $\mathbb{Q}$ sitting in some $\mathbb{Q}$-scheme $Y$ and rational functions $P_I(x,y) \in \mathbb{Q}(x,y)$ for each $I \subset T$ with the property that for almost all primes $p$*

$$(5.6) \qquad \zeta_{G,p}^*(s) = \sum_{I \subset T} c_{p,I} P_I(p, p^{-s})$$

*where*

$$c_{p,I} = \operatorname{card}\{a \in \overline{Y}(\mathbb{F}_p) : a \in \overline{E_i} \text{ if and only if } i \in I\}$$

*and $\overline{Y}$ means the reduction mod $p$ of the scheme $Y$.*

(2) *The abscissa of convergence $\alpha^*(G)$ of $\zeta_G^*(s)$ is a rational number and $\zeta_G^*(s)$ has a meromorphic continuation to $\Re(s) > \alpha^*(G) - \delta$ for some $\delta > 0$ with the property that on the line $\Re(s) = \alpha^*(G)$ the only pole is at $s = \alpha^*(G)$.*

(3) *There exist a nonnegative integer $b^*(G) \in \mathbb{N}$ and some real numbers $c, c' \in \mathbb{R}$ such that*

$$
\begin{aligned}
s_N^{*,\alpha^*(G)}(G) &= a_1^*(G) + a_2^*(G) 2^{-\alpha^*(G)} + \cdots + a_N^*(G) N^{-\alpha^*(G)} \\
&\sim c \cdot (\log N)^{b^*(G)+1} \\
s_N^*(G) &\sim c \cdot N^{\alpha^*(G)} (\log N)^{b^*(G)}
\end{aligned}
$$

*as $N \to \infty$.*

Recall that in (1) the irreducible varieties arise from the resolution of singularities of $F$, the product of all the polynomials defining the cone integral. Let us identify this $F = F_L^*$ explicitly for the cone integrals defining the zeta function of the Lie algebra $L = L(G)$ of $G$. The polynomial depends on a choice of basis for $L$ which identifies it with $\mathbb{Z}^d$. Let multiplication in $L$ be given by the bilinear mapping $\beta : \mathbb{Z}^d \times \mathbb{Z}^d \to \mathbb{Z}^d$. Let $C_j$ denote the matrix whose rows are $\mathbf{c}_i = \beta(\mathbf{e}_i, \mathbf{e}_j)$. Let $M = (m_{ij})$ be a $d \times d$ upper triangular matrix and $M^\natural = M' \operatorname{diag}(m_{22}^{-1} \cdots m_{dd}^{-1}, \ldots, m_{dd}^{-1}, 1)$ where $M'$ denotes the adjoint matrix of $M$.

*Definition* 5.8. (1) $F_L^\triangleleft(m_{ij})$ is defined to be the product of

$$m_{11}^{(d^2+1)d} \cdots m_{ii}^{(d^2+1)(d-i+1)} \cdots m_{dd}^{(d^2+1)}$$

with all the entries in the $d$ matrices $M C_j M^\natural$ for $j = 1, \ldots, d$.

(2) $F_L^{\leq}(m_{ij})$ is defined to be the product of

$$m_{11}^{(d^2+1)d} \cdots m_{ii}^{(d^2+1)(d-i+1)} \cdots m_{dd}^{(d^2+1)}$$

with all the entries in the $d$ matrices $M \left( \sum_{k=j}^{d} m_{jk} C_k \right) M^{\natural}$ for $j = 1, \ldots, d$.

The expression $m_{11}^{(d^2+1)d} \cdots m_{ii}^{(d^2+1)(d-i+1)} \cdots m_{dd}^{(d^2+1)}$ arises from the product of the two monomials in the integrand $f_0 = m_{11} \cdots m_{dd}$ and $g_0 = m_{11}^{d-1} \cdots m_{d-1 d-1}$ and the monomials $f_{ijk} = m_{11} \cdots m_{kk}$ for $i, j, k \in \{1, \ldots, d\}$.

*Definition* 5.9. For a finitely generated nilpotent group $G$ we define the polynomials $F_G^* = F_{L(G)}^*$ for $* \in \{\leq, \lhd\}$, where $L(G)$ is the associated Lie algebra of $G$.

Again there is some choice involved in the definition of the Lie algebra $L(G)$ as it is the image under log of some subgroup of finite index in $G$. However up to finitely many primes, the explicit formulas associated to the cone integrals will be the same for polynomials $F_G^* = F_{L(G)}^*$ arising from different choices of Lie algebra and basis for the Lie algebra. We refer the reader to the paper [21] of the first author and Loeser where the concept of motivic cone integrals is considered. This establishes that despite the many choices made on our way to the explicit expression (5.6), the expression is canonical. This is used in [12] to canonically associate to each nilpotent group a subring of the Grothendieck ring generated by the system of varieties $E_i$ ($i \in T$).

Although we know that the constant term of the cone integrals is $(1-p^{-1})^d$ it is instructive to see what the analysis of the resolution corresponding to $F_L^*(m_{ij})$ tells us. The constant term corresponds to when the integrand is constant which occurs when $f_0 = m_{11} \cdots m_{dd}$ is a unit, i.e. for those bases for the whole group $G$. Let $T_0 \subset T$ (where $T$ is the set indexing the irreducible components $E_i$) be defined as follows: if $i \in T_0$ then $N_i(f_0) \neq 0$. Note that if $N_i(f_0) = 0$ then $N_i(g_0) \neq 0$ since $f_0 = m_{11} \cdots m_{dd}$ and $g_0 = m_{11}^{d-1} \cdots m_{d-1 d-1}$. So those $I \subset T$ with $I \cap T_0$ nonempty are precisely those that do not contribute to the constant term. So the constant term comes from those $I \subset T \backslash T_0$. We analyse now the data in our formula associated to such $I$. If $I \subset T \backslash T_0$ then for $i \in I$, $N_i(f_0) = 0$ and since $f_{ijk} = m_{11} \cdots m_{kk}$ this also implies that $N_i(f_{ijk}) = 0$. Hence there are no conditions on our subset since the conditions are

$$\sum_{l=1}^{r} N_{i_l}(f_{ijk}) \mathrm{ord}(y_l) \leq \sum_{l=1}^{r} N_{i_l}(g_{ijk}) \mathrm{ord}(y_l).$$

This implies then that $w_I = 1$, $m_1 = m$ and $A_{k,I,1} = 0$ and $B_{k,I,1} = 1$ for $k = 1, \ldots, m$. This is the same as the case that $I = \emptyset$ which we considered at the end of Section 2. Hence the constant term is

$$p^{-m} \sum_{I \subset T \backslash T_0} c_{p,I}.$$

We know of course what the constant term should be for this integral, namely $(1 - p^{-1})^d$. Hence we get that

$$\sum_{I \subset T \backslash T_0} c_{p,I} = p^{m-d}(p-1)^d.$$

This is the same size as the number of points over $\mathbb{F}_p$ in affine space of dimension $m$ minus the hyperplanes defined by $m_{ii} = 0$. This means that the number of points has not changed outside of the union of hyperplanes $m_{ii} = 0$ after we have done the resolution. Does this mean that our variety is nonsingular outside this union of hyperplanes?

We collect together here the present knowledge we have about the possible values for $\alpha^*(G)$. Let $G$ be a torsion-free, finitely generated nilpotent group. Let $h$ denote the Hirsch length of the nilpotent group $G$ and $h_{\mathrm{ab}}$ denote the Hirsch length of the abelianisation, $h(G/G')$ .

(1) $h_{\mathrm{ab}} \leq \alpha^{\triangleleft}(G) \leq h$.

(2) $(3 - 2\sqrt{2})h - \frac{1}{2} \leq \alpha^{\leq}(G) \leq h$. The lower bound, which for large $h$ ($h \geq 17$) exceeds $h/6$, is due to Segal. The proof was extended by Klopsch to soluble groups of finite rank.

(3) If $G$ has class 2 put $m = h(Z(G))$ and $r = h(G/Z(G))$ where $Z(G)$ is the centre of $G$. Then it is proved in [24] that

$$1/2(m + r^{-1}) \leq \alpha^{\triangleleft}(G) \leq \max\{h_{\mathrm{ab}}, h(1 - r^{-1})\}.$$

(4) If $H$ has finite index in $G$ then $\alpha^*(G) = \alpha^*(H)$ (Proposition 1.8 of [24]).

The paper [24] contains a number of examples of zeta functions of class two groups calculated by Geoff Smith. There are some other examples calculated by Dermott Grenham which are recorded in [13]. In all these examples $\alpha^{\leq}(G) = \alpha^{\triangleleft}(G) = h_{\mathrm{ab}}$. However this reflects the paucity and small nature of the nilpotent groups so far considered rather than a general feature. The estimates in (3) show that for example if $G = F_2^5$ then $h_{\mathrm{ab}} = 5$ whilst $h(Z(G)) = \binom{5}{2}$. Hence

$$\alpha^{\triangleleft}(F_2^5) \geq 1/2(m + r^{-1}) > 5.$$

In [19] we shall show that the zeta function counting all subgroups in the free class two nilpotent group $F_2^3$ on three generators has abscissa of convergence at $\alpha^{\leq}(F_2^3) = 7/2$. In particular this is the first example of a group for which the abscissa of convergence is not an integer.

At present we do not know anything about the order of the poles of $\zeta_G(s)$. Our analysis above might help in giving a bound in terms of the Hirsch length since the pole at the abscissa of convergence has order given by the number of times we need to multiply by the Artin $L$-function. This is determined by the number of irreducible components of $E_I$ as $I \subset T$. The other interesting issue is some interpretation of the residue at this pole.

If the nilpotent group $G$ is not torsion-free then Proposition 1.8 of [24] can be extended to show that $\alpha^*(G) = \alpha^*(G/G^{\mathrm{tor}})$ where $G^{\mathrm{tor}}$ is the (finite) set of all torsion elements of $G$. However if we just take a finite extension of a nilpotent group then we do not know very much about the change in the abscissa of convergence. Extending even the free abelian group by a finite group can have quite a dramatic effect on the movement of the poles of $\zeta_G(s)$. As an example compare the infinite cyclic group $\mathbb{Z}$ where $\zeta_{\mathbb{Z}}(s) = \zeta(s)$ and the infinite dihedral group $D_\infty$ where $\zeta_{D_\infty}(s) = 2^{-s}\zeta(s) + \zeta(s-1)$. Extension by a finite group (here $C_2$) therefore has quite a subtle effect on the lattice of subgroups even to the extent of changing the rate of polynomial growth.

In [22] explicit examples calculated by John McDermott are recorded for the wallpaper groups which show the effect of extending $\mathbb{Z}^2$ by a finite group. It reveals how sensitive the zeta function is as the nature of the poles varies dramatically.

University of Cambridge, Cambridge, UK
*E-mail address*: dusautoy@dpmms.cam.ac.uk

Heinrich Heine Universität, Düsseldorf, Germany
*E-mail address*: fritz@math.uni-duesseldorf.de